\documentclass[preprint]{elsarticle}
\usepackage{lineno}
\modulolinenumbers[5]

\usepackage{amsmath}
\usepackage{amssymb}
\usepackage{amsthm}
\usepackage{graphicx}
\usepackage{epic,eepic,epsfig}
\usepackage{color}
\usepackage{subfigure}  
\usepackage{placeins}   
\usepackage{multirow}
\usepackage{epstopdf}
\usepackage{varwidth}
\usepackage{float}
\usepackage[table]{xcolor}

\newtheorem{theorem}{Theorem}[section]
\newtheorem{example}{Example}[section]
\newtheorem{lemma}{Lemma}[section]
\newtheorem{remark}{Remark}[section]

\newtheorem{algorithm}{Algorithm}[section]

\definecolor{tabclr}{cmyk}{0,0,1,0}

\newcommand{\bx}{\boldsymbol{x}}

\newcommand{\bff}{\boldsymbol{f}}
\newcommand{\bu}{\boldsymbol{u}}
\newcommand{\bphi}{\boldsymbol{\phi}}
\newcommand{\bV}{\boldsymbol{V}}
\newcommand{\bw}{\boldsymbol{w}}

\newcommand{\bv}{\boldsymbol{v}}
\newcommand{\bW}{\boldsymbol{W}}
\newcommand{\bU}{\boldsymbol{U}}

\newcommand{\bZ}{\boldsymbol{Z}}

\newcommand{\bg}{\boldsymbol{g}}

\newcommand{\bzeta}{\boldsymbol{\zeta}}

\newcommand{\Rmnum}[1]{\uppercase\expandafter{\romannumeral #1}} 

\begin{document}

\title{An adaptive BDDC algorithm in variational form for mortar discretizations}

\author[XTU1]{Jie Peng}
\ead{xtu\_pengjie@163.com}

\author[XTU1,XTU2]{Shi Shu}
\ead{shushi@xtu.edu.cn}

\author[XTU1]{Junxian Wang\corref{cor}}
\ead{wangjunxian@xtu.edu.cn}

\cortext[cor]{Corresponding author}
\address[XTU1]{School of Mathematics and Computational Science, Xiangtan University, Xiangtan 411105, China}
\address[XTU2]{Hunan Key Laboratory for Computation and Simulation in Science and Engineering,
Xiangtan University, Xiangtan 411105, China}

\begin{abstract}
A balancing domain decomposition by constraints (BDDC) algorithm with adaptive primal constraints in variational form
is introduced and analyzed for high-order mortar discretization of two-dimensional elliptic problems with high varying and random coefficients.
Some vector-valued auxiliary spaces and operators with essential properties are defined to describe the variational algorithm,
and the coarse space is formed by using a transformation operator on each interface.
Compared with the adaptive BDDC algorithms for conforming Galerkin approximations,
our algorithm is more simple,
because there is not any continuity constraints at subdomain vertices in the mortar method involved in this paper.
The condition number of the preconditioned system is proved to be bounded above by a user-defined tolerance and a constant which is dependent on the maximum number of interfaces per subdomain, and independent of the mesh size and the contrast of the given coefficients.
Numerical results show the robustness and efficiency of the algorithm for various model problems.
\end{abstract}


\maketitle

\section{Introduction}\label{sec:1}\setcounter{equation}{0}

Mortar methods were first introduced by Bernardi, Maday and Patera
\cite{BMP1993, BMP1994} as the discretization techniques based on
domain decomposition.
These techniques are widely applied in many scientific and
engineering computation fields, such as multi-physical models,
coupling schemes with different discretizations, problems with
non-matching grids and so on \cite{AMW1999,W2000,W2001}.
Balancing domain decomposition by constraints (BDDC) algorithms, which
were introduced by Clark R. Dohrmann \cite{Dohrmann2003}, are
variants of the balancing Neumann-Neumann algorithms for solving the
Schur complement systems. These algorithms have been extended to solve
PDE(s) discrete systems obtained by various discretization methods,
such as conforming Galerkin \cite{LW2006,BS2007,T2007}, discontinuous
Galerkin\cite{DM2007,CC2014}, and mortar methods\cite{KMW20081,KMW20082,KT2009} and so on.
However, these BDDC algorithms require a strong assumption on the coefficients in each subdomain to achieve a good performance.

To enhance the robustness, the selection of good primal constraints should be problem-dependent,
this led to adaptive algorithms for choosing primal constraints \cite{MS2007}.
Generalized eigenvalue problems with respect to the local problems per interface shared by two subdomains are used to adaptively choose primal constraints \cite{DP2012,PD2013,DSO2015,KC2015,KCW2016,KRR2016,JGC2016,PC2016,KCX2017,DV2017}.
In the work by Klawonn, Radtke and Rheinbach \cite{KRR2016}, an adaptive coarse space for the dual-primal finite element tearing and interconnecting (FETI-DP) and BDDC methods is obtained by solving generalized eigenvalue problems associate with the edge Schur complements and mass matrices.
Another class of eigenvalue problems are also introduced to construct the coarse spaces for BDDC algorithms in \cite{DP2012,KC2015},
and their eigenvalue problems are defined by using the edge Schur complements and the part of Schur complement in each subdomain.
Recently, eigenvalue problems with respect to the parallel sum (see \cite{1969AD}) have got great attention of researchers.
In Pechstein and Dohrmann \cite{PD2013}, this types of eigenvalue problems were first introduced to select the primal constraints for BDDC algorithms,
and \cite{DSO2015,KCW2016,JGC2016,KCX2017,DV2017} have extended it to elliptic problems discretized with conforming finite element methods, staggered discontinous Galerkin methods, isogeometric analysis, and vector field problems discretized with Raviart-Thomas finite elements.
However, most of the available literatures on adaptive BDDC algorithms were in algebraic form (i.e. matrices and vectors)
and adaptive BDDC algorithms for mortar discretizations have not previously been discussed in the literature.

In this paper, an adaptive BDDC algorithm in variational form for high-order mortar discretization of two dimensional elliptic problems with
high varying and random coefficients is introduced and analyzed.
Based on a vector-valued function space, we derive the Schur complement variational problem for Lagrange multiplier variable.
Then, scaling operators and transformation operators with essential properties are defined, and
a construct method of the transformation operators is presented by using the generalize eigenvalue problems with respect to the parallel sum.
Further, in contrast to the BDDC algorithms in variational form \cite{BS2007,VP2010}, by introducing some auxiliary spaces and operators, we arrive at a preconditioned adaptive BDDC algorithm in variational form for mortar discretizations.
Compared with the conforming Galerkin approximations, we emphasize that since the mortar method involved in this paper do not have any continuity constraints at subdomain vertices, this simplifies our algorithm quite a lot.
Using the characters of the involved operators, we proved that the condition number bound of the adaptive BDDC preconditioned systems is $C \Theta$,
where $C$ is a constant which depends only on the maximum number of interfaces per each subdomain, and $\Theta$ is a given tolerance.
Finally, numerical results for various model problems show the robustness of the proposed algorithms and verify the theoretical estimate in both geometrically conforming and unconforming partitions.
In particular, the algorithm with deluxe scaling matrices keeps better computational efficiency than that with multiplicity scaling matrices.

In the following, we introduce some definitions.
Assume that
$V$ and $W$ are Hilbert spaces  and $U = V \oplus W$, the operators $R:U \rightarrow
V$ and  $E:V \rightarrow U$ are separately called restriction
operator and interpolation operator refer to
\begin{align}\label{restriction-operator}
R u = v,~\forall u=v + w\in U,~\mbox{where}~v\in V, w\in W,
\end{align}
and
\begin{align}\label{interpolation-operator}
E v = v, ~\forall v \in V.
\end{align}

For a given linear operator $L$ from the Hilbert space $U$ to the Hilbert space $V$, the operator $L^T: V \rightarrow U$ is defined by
\begin{align*}
(L^T v, u) = (v, Lu),~~\forall u \in U, v\in V.
\end{align*}

The rest of this paper is organized as follows. In section 2, we present the descretization of a second order elliptic problems with mortar finite element and its corresponding Schur complement system associated with a vector-valued function space.
Some auxiliary spaces and a proper space decomposition are presented in section 3, while the adaptive BDDC algorithm is introduced in section 4.
The condition number bounds of the preconditioned system is analysed in section 5, and various numerical experiments are presented in section 6. Finally, a conclusion is given in section 7.

\section{Model problem and Schur complement system}\label{sec:2}
\setcounter{equation}{0}

Consider the following elliptic problem: find $u\in H_0^1(\Omega)$ such that
\begin{align}\label{continuous-var}
a(u, v)=(f, v),~~ \forall v\in H_0^1(\Omega),
\end{align}
where
\begin{align}\label{def-a}
a(u,v)=\int_\Omega (\rho \nabla u \cdot \nabla v + \varepsilon u v)d\bx,~~
(f,v)=\int_\Omega f vd\bx
\end{align}
and $\Omega \subset \mathbb{R}^2$ is a bounded polygonal domain, $f \in L^2(\Omega)$, the bounded coefficients $\varepsilon \ge 0$ and $\rho \ge \rho_{\min} >0$
can be random and has  high contrast  in $\Omega$.

In the following, we define a mortar discrete problem of \eqref{continuous-var} based on a nonoverlapping domain decompositions.

We decompose the given region $\Omega$ into polyhedral subdomains $\Omega_i(i=1,\cdots,N)$, which satisfy
\begin{align*}
\bar{\Omega}=\bigcup\limits_{k=1}^{N} \bar{\Omega}_k ~~\mbox{with}~~\Omega_i \cap \Omega_j = \emptyset,~~i\neq j,
\end{align*}
$\partial\Omega_i\cap\partial\Omega_j(i \neq j)$ is either empty, a vertex or a common edge,
and let $d_i$ be the diameter of $\Omega_i$.
Each subdomain $\Omega_i$ is associated with a regular or quasi-uniform triangulation $\mathcal{T}_i$,
where the mesh size of $\mathcal{T}_i$ is denoted by $h_i$.

Denote the $\mathcal{P}_s$ Lagrange finite element space associated with $\mathcal{T}_i$ by
\begin{align*}
X(\Omega_i) = \{v \in C(\overline{\Omega}_i): v \big{|}_\tau \in \mathcal{P}_s, \forall \tau  \in \mathcal{T}_i , ~v \big{|}_{\partial \Omega_i \cap \partial \Omega} =0 \},
\end{align*}
where $\mathcal{P}_s$ denotes the set of all polynomials
of degree less than or equal to $s$, and $s \in \mathbb{Z}^+$.

If $\partial \Omega_i \cap \partial \Omega_j(i\neq j)$ is a common edge, we call it an interface.
Each interface $\partial \Omega_i \cap \partial \Omega_j(i\neq j)$ is associated with a one-dimensional triangulation, provided either from $\mathcal{T}_i$ or $\mathcal{T}_j$. Since the triangulations in any different subdomains are independent of each other, they are generally do not match at the interfaces. For convenience, we denote the interface by $\Gamma_{ij}$ and $\Gamma_{ji}$, respectively, when the triangulation is given by $\mathcal{T}_i$ and $\mathcal{T}_j$.
Further, let $\Gamma=\cup\Gamma_{ij}$, $M$ denote the number of interfaces.

Denote
\begin{align*}
\mathcal{E} := \{\Gamma_{ij}:~1 \le i \le N, j \in \mathcal{E}_i\},~~\mathcal{M} :=\{\Gamma_{ij} \in \mathcal{E}:~h_i \ge h_j\},
\end{align*}
where
$
\mathcal{E}_i := \{j:~\partial \Omega_i \cap \partial \Omega_j ~\mbox{is a common edge},~\mbox{for}~1\le j\le N~\mbox{and}~j\neq i\}.
$


Let the interfaces in $\mathcal{M}$ denote the nonmortars, and those of $\mathcal{E}\backslash
\mathcal{M}$ the mortars (see \cite{W2000}). The discrete Lagrange multiplier space will be associated with the nonmortars.
Since there is one-to-one correspondence between
element $\Gamma_{ij}$ in $\mathcal{M}$ and the interface,
we can denote $\mathcal{M}$ by $\{F_k: k=1,\cdots,M\}$,
where $F_k$ is called the interface with global index $k$. For any given subdomain $\Omega_i$, let
\begin{align}\label{def-M-i}
\mathcal{M}_i := \{k:~F_k \subset \partial \Omega_i \backslash \partial \Omega,~\mbox{for}~1\le k \le M\}.
\end{align}
%
%
%
%
%


Let $M(F_k)$ denote the standard Lagrange multiplier space with respect to the nonmortar edge $F_k$ (see \cite{KMW20081,BMP1993,BMP1994}) and $n_k = dim(M(F_k))$.
%
%
%
Define the extension space of $X(\Omega_i)$ and $M(F_k)$ by
\begin{align*}
X^{(i)} = E^{(i)}(X(\Omega_i))~\mbox{and}~M_{F_k} = E_k(M(F_k)),
\end{align*}
where the trivial extension operator $E^{(i)}$ and $E_k$ satisfy
\begin{align*}
E^{(i)} v = \left\{\begin{array}{ll}
v  & \mbox{on}~\bar{\Omega}_i\\
0  & \mbox{on}~\Omega \backslash \bar{\Omega}_i
\end{array}\right.,~\forall v\in X(\Omega_i)~\mbox{and}~
E_k \lambda = \left\{\begin{array}{ll}
\lambda  &\mbox{on}~\bar{F}_k\\
0  & \mbox{on}~\Gamma \backslash \bar{F}_k
\end{array}\right.,~\forall \lambda \in M(F_k).
\end{align*}


Denote the direct sum of $X^{(i)}(i=1,\cdots,N)$ and
$M_{F_k}(k=1,\cdots,M)$ respectively as
\begin{align*}
X_h = \oplus_{i=1}^N X^{(i)}~\mbox{and}~M_h = \oplus_{k = 1}^M M_{F_k}.
\end{align*}

Let $\bV_h$ denotes the vector-valued function space $X_h \times M_h$.
The mortar finite element approximation of problem \eqref{continuous-var} is as follows(the case for $\mathcal{P}_1$ Lagrange finite element space see \cite{W2000}).

Find $(u, \lambda) \in \bV_h$ such that
\begin{align}\label{continue-variable-from}
\left\{
\begin{array}{lcll}
 \tilde{a}(u,v)  + b(v,\lambda) &=& (f,v), & \forall v\in X_h,\\
 b(u, q)            &=& 0, & \forall q \in M_h,
\end{array}\right.
\end{align}
where
\begin{align}\label{def-AUV-000}
\tilde{a}(u,v) = \sum\limits_{i=1}^N \tilde{a}_i(u,v),~~ \tilde{a}_i(u,v) = \int _{\Omega_i} (\rho \nabla u \cdot
\nabla v + \varepsilon u v)d\bx,
\end{align}
and
\begin{align}\label{def-sigma}
b(v, \mu) = \sum\limits_{i=1}^N b_i(v, \mu),~~b_i(v, \mu) = \sum \limits_{\Gamma_{ir}
\subset \partial \Omega_i \backslash \partial \Omega} \int_{\Gamma_{ir}}
(\sigma_{ir} v|_{\Omega_i}  \mu ) dS,~~
\sigma_{ir} = \left\{\begin{array}{ll}
1, & i < r,\\
-1, & i > r.
\end{array}\right.
\end{align}

\begin{remark}
When $\Omega_i$ is an internal subdomain and $\varepsilon =0$, the bilinear form $\tilde{a}_i(\cdot, \cdot)$ which is symmetric positive semi-definite
can be regularized and transformed to a symmetric positive definite (SPD) form  (see \cite{HU2007}).
Therefore, we always assume that $\tilde{a}_i(\cdot,\cdot)(i=1,\cdots, N)$ are coercive on $X_h$.
\end{remark}

The Schur complement of the system \eqref{continue-variable-from} is
\begin{align}\label{Schur-complement-pro}
S \lambda = g,
\end{align}
where
\begin{align*}
S = \bar{B} \bar{A}^{-1} \bar{B}^T~\mbox{and}~g = \bar{B} \bar{A}^{-1} f,
\end{align*}
here the operator $\bar{A}: X_h \rightarrow X_h$ and  $\bar{B}: X_h \rightarrow M_h$ such that
\begin{align*}
(\bar{A} u, v) = \tilde{a}(u,v),~\forall u, v \in X_h~\mbox{and}~
(\bar{B}v, \lambda) = b(v,\lambda),~\forall v\in X_h, \lambda \in
M_h.
\end{align*}
In order to discuss the adaptive BDDC preconditioner in variational form for solving the mortar discretizations \eqref{continue-variable-from} restricted to the coupling of $\mathcal{P}_s$-Lagrangian finite elements, we need to derive the corresponding variational problem of
the scalar Schur complement system \eqref{Schur-complement-pro} on a vector-valued function space.


For any $\bu = (u,\lambda), \bv = (v, q)\in \bV_h$, we introduce a bilinear form
\begin{align}\label{def-AUV}
A(\bu, \bv) = \sum\limits_{i=1}^N A_i(\bu, \bv),~\mbox{where}~A_i(\bu, \bv) =
\tilde{a}_i(u,v) + b_i(v,\lambda) + b_i(u,q)~\mbox{for}~1\le i \le N,
\end{align}
%
%
%
then the saddle point problem \eqref{continue-variable-from} is equivalent to the following variational problem: find $\bu \in \bV_h$ such that
\begin{align}\label{continue-variable-from-2}
A(\bu, \bv) = (\bff, \bv),~\forall~\bv \in \bV_h,
\end{align}
where
\begin{align}\label{def-AUV-i}
\bff = (f, 0),~(\bff, \bv) = (f, v), ~\forall~\bv = (v, q)\in
\bV_h.
\end{align}




Define the vector-valued function spaces
\begin{align}\label{def-V-1}
\bV_{I} = \oplus_{i=1}^N \bV_I^{(i)},~~\bV_I^{(i)} = \{(v,0):~v\in X^{(i)}\},~i=1,\cdots,N,
\end{align}
and
\begin{align}\label{def-V-2}
\bV_{F_k}=\{(0,\lambda):~\lambda \in M_{F_k}\}, ~k=1,\cdots M.
\end{align}


For any given subdomain $\Omega_i$, let
\begin{align}\label{def-V-3}
\bV^{(i)} = (\oplus_{k \in \mathcal{M}_i} \bV_{F_k}) \oplus \bV_I^{(i)},
\end{align}
where $\mathcal{M}_i$ is defined in \eqref{def-M-i}.


For any interface $F_k$, by using the multiplier basis functions $\{\varphi^{k}_l\}_{l=1}^{n_k}$ of $M_{F_k}$,
we can define a function vector  
\begin{align}\label{def-Phik}
\Phi^k = (\bphi^{k}_1,\cdots,\bphi^{k}_{n_k})^T,
\end{align}
where the $l$-th vector-valued function $\bphi^{k}_l = (\phi^{k}_{l},
\psi^{k}_{l}) \in \bV_I \oplus \bV_{F_k}$ satisfies that
\begin{align}\label{def-Hw}
\left\{\begin{array}{l}
A(\bphi^{k}_l, \bv) = 0, ~~\forall \bv \in \bV_I,\\
\psi^{k}_{l}|_{F_k} = \varphi^{k}_l.
\end{array}\right.
\end{align}

Utilizing the vector $\Phi^k$ defined in \eqref{def-Phik}, we can define vector-valued function spaces
\begin{align}\label{def-W}
\hat{\bW} = \oplus_{k=1}^M \bW_{k},~\mbox{where}~\bW_{k} = span\{\bphi^{k}_1,\cdots,\bphi^{k}_{n_k}\},
\end{align}
and the variational form of the Schur complement system for \eqref{continue-variable-from-2} can be expressed as:
find $\hat{\bw} = (w, \lambda) \in \hat{\bW}$ such that
\begin{align}\label{def-schur-bilinear-form-trans-old}
A(\hat{\bw}, \hat{\bv}) = (\bff, \hat{\bv}),~~\forall \hat{\bv} \in \hat{\bW}.
\end{align}


Obviously, the second component of the solution to the above variational problem, i.e. $\lambda$, is also the solution of the Schur complement system \eqref{Schur-complement-pro}.


Let $\hat{S}: \hat{\bW} \rightarrow \hat{\bW}$ be the Schur complement operator defined by
\begin{align}\label{def-operator-hat-S}
(\hat{S} \hat{\bu}, \hat{\bv}) = A(\hat{\bu}, \hat{\bv}),~~\forall \hat{\bu}, \hat{\bv} \in \hat{\bW}.
\end{align}
We can rewrite \eqref{def-schur-bilinear-form-trans-old} as 
\begin{align}\label{def-schur-system}
\hat{S} \hat{\bw}  = Q_{\hat{W}} \bff,
\end{align}
where $Q_{\hat{W}}:(L^2(\Omega),L^2(\Gamma)) \rightarrow \hat{\bW}$ is the $L^2$ projection operator.


In order to give an adaptive BDDC preconditioner for solving the Schur complement system \eqref{def-schur-system},
some auxiliary spaces and a proper decomposition of $\hat{\bW}$ are presented in the next section.


\section{Some auxiliary spaces and space decomposition}\label{sec:3}\setcounter{equation}{0} 

For any given interface $F_k$ $(k=1,\cdots,M)$, we construct a new set of basis functions of the space $\bW_{k}$ defined in \eqref{def-W}.


We always assume that $F_k$ be the interface shared by $\Omega_i$ and $\Omega_j$.
For $\nu=i,j$, let $\bV^{(\nu)}_I, \bV_{F_k}$ and $\bV^{(\nu)}$ be the spaces defined in \eqref{def-V-1},
\eqref{def-V-2} and \eqref{def-V-3} respectively.
By using the basis functions $\{\varphi^{k}_l\}_{l=1}^{n_k}$, the vectors
\begin{align}\label{def-Phi-barPhi-k-nu}
\Phi^{k,\nu} = (\bphi^{k,\nu}_1, \cdots, \bphi^{k,\nu}_{n_k})^T~\mbox{and}~\bar{\Phi}^{k,\nu} = (\bar{\bphi}^{k,\nu}_1, \cdots, \bar{\bphi}^{k,\nu}_{n_k})^T
\end{align}
can be defined similarly to $\Phi^k$ in \eqref{def-Phik}, where
$\bphi^{k,\nu}_l = (\phi^{k,\nu}_{l}, \psi^{k,\nu}_{l}) \in \bV_I^{(\nu)} \oplus \bV_{F_k}$,
  $\bar{\bphi}^{k,\nu}_l = (\bar{\phi}^{k,\nu}_{l}, \bar{\psi}^{k,\nu}_{l})\in \bV^{(\nu)}$ satisfy that
\begin{align}\label{def-Hw-i}
\left\{\begin{array}{l}
A_{\nu}(\bphi^{k,\nu}_l, \bv) = 0, ~~\forall \bv \in \bV^{(\nu)}_I\\
\psi^{k,\nu}_{l}|_{F_k} = \varphi^{k}_l
\end{array}\right.,~~l=1,\cdots,n_k,
\end{align}
and
\begin{align}\label{def-barHw-i-Fim}
\left\{\begin{array}{l}
A_{\nu}(\bar{\bphi}^{k,\nu}_l, \bv) = 0, ~~\forall \bv \in \bV^{(\nu)} \backslash \bV_{F_k}\\
\bar{\psi}^{k,\nu}_{l}|_{F_k} = \varphi^{k}_l
\end{array}
\right.,~~l=1,\cdots,n_k.
\end{align}


Define the auxiliary vector-valued function spaces ($\nu=i,j$)
\begin{align}\label{def-W-k-nu}
& \bW_{k}^{(\nu)}=span\{\bphi^{k,\nu}_1,\cdots, \bphi^{k,\nu}_{n_k}\},~\bar{\bW}_{k}^{(\nu)} = span\{\bar{\bphi}^{k,\nu}_1,\cdots, \bar{\bphi}^{k,\nu}_{n_k}\},\\\label{def-Z-k-nu}
& \bZ_{k}^{(\nu)} = span\{\bar{\bphi}^{k,\nu}_1 - \bphi^{k,\nu}_1,\cdots, \bar{\bphi}^{k,\nu}_{n_k} - \bphi^{k,\nu}_{n_k}\}.
\end{align}
Using \eqref{def-Phi-barPhi-k-nu}, \eqref{def-Hw-i} and \eqref{def-barHw-i-Fim}, we derive
\begin{align}\label{def-w1-property-1-000}
A_{\nu}(\bw, \bv) = 0,~~\forall \bw \in \bar{\bW}_{k}^{(\nu)},~\bv \in \bU,~\nu=i,j,
\end{align}
where $\bU = \bZ_{k}^{(\nu)}~\mbox{or}~\bW_{m}^{(\nu)}$, $m \in \mathcal{M}_{\nu}$ and $m \neq k$.


Let $D_{F_{k}}^{(\nu)}: \bU \rightarrow \bU$ $(\nu=i,j)$ be the scaling operator, where $\bU = \bW_{k}^{(i)}~\mbox{or}~\bW_{k}^{(j)}$,
and satisfy that for all $\bw = \vec{w}^T \Psi $ with $\vec{w} \in \mathbb{R}^{n_k}$ and $\Psi = \Phi^{k,i}~\mbox{or}~\Phi^{k,j}$, we have
\begin{align}\label{def-D-operator-matrix}
D_{F_{k}}^{(\nu)} \bw = \vec{w}^T (\vec{D}_{F_{k}}^{(\nu)})^T \Psi,
\end{align}
where the $n_k\times n_k$ scaling matrix $\vec{D}_{F_{k}}^{(\nu)}$ is nonsingular, and
\begin{align}\label{def-D-trans}
D_{F_{k}}^{(i)} + D_{F_{k}}^{(j)}  = I,~\mbox{where $I$ is the identity operator.}
\end{align}

Two of the most frequently used formulas of the scaling matrices $\vec{D}_{F_{k}}^{(\nu)}(\nu=i,j)$ are (see \cite{RF1999,DP2012})
\begin{align}\label{DiFD}
\vec{D}^{(i)}_{F_k}=\frac{1}{2}\vec{I},~\vec{D}^{(j)}_{F_k}=\frac{1}{2}\vec{I},
\end{align}
and
\begin{align}\label{DiFD-deluxe}
\vec{D}^{(i)}_{F_k}=(\vec{S}^{(i)}_{F_k} +
\vec{S}^{(j)}_{F_k})^{-1}\vec{S}^{(i)}_{F_k},~
\vec{D}^{(j)}_{F_k}=(\vec{S}^{(i)}_{F_k} +
\vec{S}^{(j)}_{F_k})^{-1}\vec{S}^{(j)}_{F_k},
\end{align}
where $\vec{I}$ denotes the $n_k \times n_k$ identity matrix, and
\begin{align}\label{Sij-def}
\vec{S}^{(\nu)}_{F_k}  = (a_{l,m}^{(\nu)})_{n_k\times n_k},~a_{l,m}^{(\nu)} = A_{\nu}(\bphi^{k,\nu}_m, \bphi^{k,\nu}_l),~l,m=1,\cdots, n_k,~\nu=i,j.
\end{align}

The matrices defined in \eqref{DiFD} and \eqref{DiFD-deluxe} are usually called \emph{multiplicity scaling matrices} and \emph{deluxe
scaling matrices}, respectively.

For any given positive real number $\Theta \ge 1$,
using the scaling operator $D_{F_{k}}^{(\nu)}$ $(\nu=i,j)$
and the function spaces defined in \eqref{def-W} and \eqref{def-W-k-nu},
we can define a linear transformation operator $T_{F_k}:\bU \rightarrow \bU$ ($\bU = \bW_{k},
\bW_{k}^{(\nu)}, \bar{\bW}_{k}^{(\nu)}, \nu=i,j$) such that for each $\bw =
\vec{w}^T \Psi$ with $\vec{w} \in \mathbb{R}^{n_k}$ and $\Psi = \Phi^k,
\Phi^{k,\nu}$ or $\bar{\Phi}^{k,\nu}$ $(\nu=i,j)$, we have
\begin{align}\label{def-bar-Phi-k-Delta1-nu}
T_{F_k} \bw  = \vec{w}^T (\vec{T}_{F_k})^T \Psi,
\end{align}
where the $n_k \times n_k$ transformation matrix $\vec{T}_{F_k}$ is nonsingular,
and satisfies the following condition: for any given $\vec{w}_{\Delta} =
(w_1,\cdots,w_{n_{\Delta}^k})^T \in \mathbb{R}^{n_{\Delta}^k},
\vec{w}_{\Pi} = (v_1,\cdots,v_{n_{\Pi}^k})^T \in
\mathbb{R}^{n_{\Pi}^k}$, $n_k = n_{\Delta}^k + n_{\Pi}^k$, we have
\begin{align}\label{def-TFkDelta-2-equ-1}
A_i(D_{F_k}^{(j)} \bw_{k,\Delta}^{(i)},D_{F_k}^{(j)}\bw_{k,\Delta}^{(i)}) +  A_j(D_{F_k}^{(i)} \bw_{k,\Delta}^{(j)}, D_{F_k}^{(i)} \bw_{k,\Delta}^{(j)})
\le \Theta A_i(\bar{\bw}_{k,\Delta}^{(i)} + \bar{\bw}_{k,\Pi}^{(i)}, \bar{\bw}_{k,\Delta}^{(i)} + \bar{\bw}_{k,\Pi}^{(i)}),
\end{align}
here
\begin{align}\label{def-TFkDelta-2-equ-2}
\bar{\bw}_{k,\Delta}^{(i)} = \sum\limits_{l=1}^{n_{\Delta}^k} w_l T_{F_k} \bar{\bphi}_l^{k,i},~\bar{\bw}_{k,\Pi}^{(i)} = \sum\limits_{l=1}^{n_{\Pi}^k} v_l T_{F_k} \bar{\bphi}_{n_{\Delta}^k+l}^{k,i},~\bw_{k,\Delta}^{(\nu)} = \sum\limits_{l=1}^{n_{\Delta}^k} w_l T_{F_k} \bphi_l^{k,\nu},~\nu=i,j.
\end{align}

We now give a way to construct the linear operator $T_{F_k}$.
Using the bilinear form $A_{\nu}(\cdot,\cdot)$ defined in \eqref{def-AUV}, and the basis functions $\{\bar{\bphi}^{k,\nu}_l\}_{l=1}^{n_k}$ defined in \eqref{def-Phi-barPhi-k-nu}, we can define two matrices via
\begin{align}\label{Sij-def-2}
\vec{\bar{S}}^{(\nu)}_{F_k}  = (b_{l,m}^{(\nu)})_{n_k\times n_k},~b_{l,m}^{(\nu)} = A_{\nu}(\bar{\bphi}^{k,\nu}_m, \bar{\bphi}^{k,\nu}_l),~l,m=1,\cdots, n_k,~\nu=i,j,
\end{align}
and their parallel sum  (see \cite{1969AD})
$$\vec{\bar{S}}_{F_{k}}^{(i)}:\vec{\bar{S}}_{F_{k}}^{(j)} =
\vec{\bar{S}}_{F_{k}}^{(j)}(\vec{\bar{S}}_{F_{k}}^{(i)}+\vec{\bar{S}}_{F_{k}}^{(j)})^{\dagger}\vec{\bar{S}}_{F_{k}}^{(i)},$$
where
$(\vec{\bar{S}}_{F_{k}}^{(i)}+\vec{\bar{S}}_{F_{k}}^{(j)})^{\dagger}$
is a pseudo inverse of the matrix
$\vec{\bar{S}}_{F_{k}}^{(i)}+\vec{\bar{S}}_{F_{k}}^{(j)}$.


Since $\vec{\bar{S}}_{F_{k}}^{(\nu)}(\nu=i,j)$ are both SPD,  $\vec{\bar{S}}_{F_{k}}^{(i)}:\vec{\bar{S}}_{F_{k}}^{(j)}$
is also SPD and satisfies
\begin{align}\label{parallel-property}
\vec{\bar{S}}_{F_{k}}^{(i)}:\vec{\bar{S}}_{F_{k}}^{(j)} \le \vec{\bar{S}}_{F_{k}}^{(\nu)},~\nu=i,j.
\end{align}


Introducing a generalized eigenvalue problem (see \cite{PD2013,KCW2016,PC2016,KCX2017})
\begin{align}\label{eig-pro-intro}
((\vec{D}_{F_k}^{(i)})^T \vec{S}_{F_k}^{(j)} \vec{D}_{F_k}^{(i)} + (\vec{D}_{F_k}^{(j)})^T \vec{S}_{F_k}^{(i)} \vec{D}_{F_k}^{(j)}) \vec{v}
 = \lambda \vec{\bar{S}}_{F_{k}}^{(i)}:\vec{\bar{S}}_{F_{k}}^{(j)}\vec{v},
\end{align}
where $\vec{v} \in \mathbb{R}^{n_k}$,
$\vec{D}_{F_k}^{(\nu)}(\nu=i,j)$ are the scaling matrices and
$\vec{S}^{(\nu)}_{F_k}(\nu=i,j)$ are defined in \eqref{Sij-def}.


Let
\begin{align*}
\lambda_1 \le \lambda_2 \le \cdots \le \lambda_{n_{\Delta}^k} \le \Theta \le \lambda_{n_{\Delta}^k+1} \le \cdots \le \lambda_{n_k}
\end{align*}
be the eigenvalues of \eqref{eig-pro-intro}, where $n_{\Delta}^k$ is a non-negative integer, and $\Theta \ge 1$ is a given tolerance in \eqref{def-TFkDelta-2-equ-1}.


Denote the $n_k\times n_k$ transformation matrix
\begin{align}\label{def-T-1}
\vec{T}_{F_k} = (\vec{T}_{\Delta}^{F_k}, \vec{T}_{\Pi}^{F_k}),
\end{align}
where
\begin{align}\label{def-T-Delta1-Pi1}
\vec{T}_{\Delta}^{F_k}:= (\vec{v}_1,\cdots,\vec{v}_{n_{\Delta}^k}),~~\vec{T}_{\Pi}^{F_k}:= (\vec{v}_{n_{\Delta}^k+1},\cdots,\vec{v}_{n_k}),
\end{align}
here $\vec{v}_l(l=1,\cdots,n_k)$ are the generalized eigenvectors  of \eqref{eig-pro-intro} corresponding to $\lambda_l$.


Using the above matrix $\vec{T}_{F_k}$, we can obtain the operator $T_{F_k}$ defined in \eqref{def-bar-Phi-k-Delta1-nu}.
Next, we want to verify that it satisfies \eqref{def-TFkDelta-2-equ-1}.


 For the special choice of $\bw=\bphi^{k,\nu}_l$$(l=1,\cdots,n_{\Delta}^k)$ in \eqref{def-bar-Phi-k-Delta1-nu}, it is easy to know that  $\vec{w} = (\delta_{l,1},\cdots,\delta_{l,n_k})^T$, where $\delta_{l,m}(m=1,\cdots,n_k)$ are the Kronecker delta.
From this and utilizing \eqref{def-T-1} and
\eqref{def-T-Delta1-Pi1}, it follows that
\begin{align}\label{basis-function-1}
(T_{F_k} \bphi^{k,\nu}_1,\cdots,T_{F_k} \bphi^{k,\nu}_{n_{\Delta}^k})^T = (\vec{T}_{\Delta}^{F_k})^T \Phi^{k,\nu},~\nu=i,j.
\end{align}
Similarly,
\begin{align}\label{basis-function-2}
(T_{F_k} \bar{\bphi}^{k,i}_{1},\cdots,T_{F_k} \bar{\bphi}^{k,i}_{n_{\Delta}^k})^T = (\vec{T}_{\Delta}^{F_k})^T \bar{\Phi}^{k,i}, (T_{F_k} \bar{\bphi}^{k,i}_{n_{\Delta}^k+1},\cdots,T_{F_k} \bar{\bphi}^{k,i}_{n_k})^T = (\vec{T}_{\Pi}^{F_k})^T \bar{\Phi}^{k,i}.
\end{align}
%
%
Then, 
we can rewrite the functions in \eqref{def-TFkDelta-2-equ-2} as
\begin{align}\label{wdelta-wpi-def}
\bar{\bw}_{k,\Delta}^{(i)} = \vec{w}_{\Delta}^T (\vec{T}_{\Delta}^{F_k})^T \bar{\Phi}^{k,i},
~\bar{\bw}_{k,\Pi}^{(i)} = \vec{w}_{\Pi}^T (\vec{T}_{\Pi}^{F_k})^T \bar{\Phi}^{k,i},
~\bw_{k,\Delta}^{(\nu)} &= \vec{w}_{\Delta}^T (\vec{T}_{\Delta}^{F_k})^T \Phi^{k,\nu},~\nu=i,j.
\end{align}


 By using \eqref{wdelta-wpi-def}, \eqref{def-D-operator-matrix}, \eqref{def-T-Delta1-Pi1}, \eqref{eig-pro-intro} and \eqref{parallel-property},
 and note that the eigenvectors $\vec{v}_l(l=1,\cdots,n_{\Delta}^k)$ are orthogonality and their corresponding eigenvalue $\lambda_l \le \Theta$, $\vec{\bar{S}}_{F_{k}}^{(i)}:\vec{\bar{S}}_{F_{k}}^{(j)}$ is SPD, we have
\begin{align*}
& A_i(D_{F_k}^{(j)} \bw_{k,\Delta}^{(i)},D_{F_k}^{(j)}\bw_{k,\Delta}^{(i)}) +  A_j(D_{F_k}^{(i)} \bw_{k,\Delta}^{(j)}, D_{F_k}^{(i)} \bw_{k,\Delta}^{(j)})\\\nonumber
&= A_i(\vec{w}_{\Delta}^T (\vec{T}_{\Delta}^{F_k})^T (\vec{D}_{F_k}^{(j)})^T \Phi^{k,i},\vec{w}_{\Delta}^T (\vec{T}_{\Delta}^{F_k})^T (\vec{D}_{F_k}^{(j)})^T \Phi^{k,i})
+  A_j(\vec{w}_{\Delta}^T (\vec{T}_{\Delta}^{F_k})^T (\vec{D}_{F_k}^{(i)})^T \Phi^{k,j}, \vec{w}_{\Delta}^T (\vec{T}_{\Delta}^{F_k})^T (\vec{D}_{F_k}^{(i)})^T \Phi^{k,j})\\\nonumber
&= \vec{w}_{\Delta}^T (\vec{T}_{\Delta}^{F_k})^T \left((\vec{D}_{F_k}^{(j)})^T \vec{S}_{F_k}^{(i)} \vec{D}_{F_k}^{(j)} + (\vec{D}_{F_k}^{(i)})^T \vec{S}_{F_k}^{(j)}\vec{D}_{F_k}^{(i)}\right) \vec{T}_{\Delta}^{F_k} \vec{w}_{\Delta} \\\nonumber
&= \Theta \vec{w}_{\Delta}^T (\vec{T}_{\Delta}^{F_k})^T (\vec{\bar{S}}_{F_{k}}^{(i)}:\vec{\bar{S}}_{F_{k}}^{(j)}) \vec{T}_{\Delta}^{F_k} \vec{w}_{\Delta}\\
&= \Theta (\vec{T}_{\Delta}^{F_{k}}\vec{w}_{\Delta} + \vec{T}_{\Pi}^{F_{k}}\vec{w}_{\Pi})^T (\vec{\bar{S}}_{F_{k}}^{(i)}:\vec{\bar{S}}_{F_{k}}^{(j)}) (\vec{T}_{\Delta}^{F_{k}}\vec{w}_{\Delta} + \vec{T}_{\Pi}^{F_{k}}\vec{w}_{\Pi})\\\nonumber
&\le \Theta (\vec{T}_{\Delta}^{F_{k}}\vec{w}_{\Delta} + \vec{T}_{\Pi}^{F_{k}}\vec{w}_{\Pi})^T \vec{\bar{S}}_{F_{k}}^{(i)} (\vec{T}_{\Delta}^{F_{k}}\vec{w}_{\Delta} + \vec{T}_{\Pi}^{F_{k}}\vec{w}_{\Pi})\\\nonumber
&= \Theta
A_i(\bar{\bw}_{k,\Delta}^{(i)} + \bar{\bw}_{k,\Pi}^{(i)}, \bar{\bw}_{k,\Delta}^{(i)} + \bar{\bw}_{k,\Pi}^{(i)})
\end{align*}
Then \eqref{def-TFkDelta-2-equ-1} holds.
This completes the construction of the operator $T_{F_k}$.

Using the linear operator $T_{F_k}$ defined in \eqref{def-bar-Phi-k-Delta1-nu}, we can obtain a new set of basis functions of $\bW_{k}$ as follows
\begin{align}\label{basis-1}
\{\hat{\bphi}_l^k:=T_{F_k} \bphi_l^k: 1 \le l \le n_k\}.
\end{align}


Based on the basis functions described above, we can decompose the space $\bW_{k}$ into
\begin{align}\label{def-hatWk}
\bW_{k} = \bW_{k,\Delta}\oplus  \bW_{k,\Pi},
\end{align}
where
\begin{align}\label{def-W-k}
\bW_{k,\Delta} = span\{\hat{\bphi}^{k}_1, \cdots, \hat{\bphi}^{k}_{n_{\Delta}^k}\}~\mbox{and}~
\bW_{k,\Pi} = span\{\hat{\bphi}^{k}_{n_{\Delta}^k+1}, \cdots, \hat{\bphi}^{k}_{n_k}\}.
\end{align}


Then using \eqref{def-W} and \eqref{def-hatWk}, a decomposition of the space $\hat{\bW}$ can be obtained as follows
\begin{align}\label{def-hatW}
\hat{\bW} = \oplus_{k=1}^M \bW_{k,\Delta} \oplus \bW_{\Pi},
\end{align}
where the coarse-level and primal variable space
\begin{align}\label{def-W-Pi}
\bW_{\Pi} = \oplus_{k=1}^M \bW_{k,\Pi}.
\end{align}


Similarly, using the linear operator $T_{F_k}$, we can get a new set of basis functions separately for the auxiliary spaces $\bW_{k}^{(\nu)}$ and
$\bar{\bW}_{k}^{(\nu)}$ $(\nu=i,j)$ defined in \eqref{def-W-k-nu} as
\begin{align}\label{basis-2}
\{\tilde{\bphi}_l^{k,\nu}:=T_{F_k} \bphi_l^{k,\nu}: 1 \le l \le n_k\}~\mbox{and}~\{\tilde{\bar{\bphi}}_l^{k,\nu}:=T_{F_k} \bar{\bphi}_l^{k,\nu}: 1 \le l \le n_k\}.
\end{align}
Then, we decompose $\bW_{k}^{(\nu)}$ and $\bar{\bW}_{k}^{(\nu)}$ into
\begin{align}\label{def-hatWk-nu}
\bW_{k}^{(\nu)} = \bW_{k,\Delta}^{(\nu)} \oplus  \bW_{k,\Pi}^{(\nu)}~\mbox{and}~\bar{\bW}_{k}^{(\nu)} = \bar{\bW}_{k,\Delta}^{(\nu)} \oplus  \bar{\bW}_{k,\Pi}^{(\nu)},
\end{align}
where
\begin{align}\label{def-W-k-Delta-nu}
&\bW_{k,\Delta}^{(\nu)} = span\{\tilde{\bphi}^{k,\nu}_1, \cdots, \tilde{\bphi}^{k,\nu}_{n_{\Delta}^k}\},
~\bW_{k,\Pi}^{(\nu)}    = span\{\tilde{\bphi}^{k,\nu}_{n_{\Delta}^k+1}, \cdots, \tilde{\bphi}^{k,\nu}_{n_k}\},\\\label{def-barW-k-nu}
&\bar{\bW}_{k,\Delta}^{(\nu)} = span\{\tilde{\bar{\bphi}}^{k,\nu}_1, \cdots, \tilde{\bar{\bphi}}^{k,\nu}_{n_{\Delta}^k}\},
~\bar{\bW}_{k,\Pi}^{(\nu)} = span\{\tilde{\bar{\bphi}}^{k,\nu}_{n_{\Delta}^k+1}, \cdots, \tilde{\bar{\bphi}}^{k,\nu}_{n_k}\}.
\end{align}


By using $\bW_{\Pi}$ defined in \eqref{def-W-Pi}, and $\bW_{k,\Delta}^{(i)}$ $(k\in \mathcal{M}_i, 1\le i\le N)$ defined in \eqref{def-W-k-Delta-nu}, we can arrive at an another auxiliary space
\begin{align}\label{def-tildeW}
\tilde{\bW} = \tilde{\bW}_{\Delta} \oplus \bW_{\Pi},
\end{align}
where
\begin{align}\label{def-tildeW-Delta}
\tilde{\bW}_{\Delta} = \oplus_{i=1}^N \bW_{\Delta}^{(i)},
~\bW_{\Delta}^{(i)} = \oplus_{k \in \mathcal{M}_i} \bW_{k,\Delta}^{(i)},~i=1,\cdots,N.
\end{align}


In the next section, we will present the adaptive BDDC preconditioner for solving the Schur complement system \eqref{def-schur-system}
by using the decomposition \eqref{def-hatW} of space $\hat{\bW}$ and some auxiliary spaces introduced in this section.


\section{Adaptive BDDC preconditioner}\label{sec:4}\setcounter{equation}{0}

For any given subdomain $\Omega_i (i=1,\cdots,N)$ and its interface $F_k$ $(k \in \mathcal{M}_i)$,
let $T_k^{(i)}: \bW_k \rightarrow \bW_k^{(i)}$ be the linear basis transformation operator such that
\begin{align}\label{def-T-k-i}
T_k^{(i)} \hat{\bphi}^k_l = \tilde{\bphi}^{k,i}_l,~l=1,\cdots, n_k,
\end{align}
where the basis functions
$\{\hat{\bphi}^k_l\}_{l=1}^{n_k}$ and
$\{\tilde{\bphi}^{k,i}_l\}_{l=1}^{n_k}$ are separately defined in
\eqref{basis-1} and \eqref{basis-2}.

For any given $\tilde{\bu}, \tilde{\bv} \in \tilde{\bW}$, by the definition \eqref{def-tildeW}, we  have
\begin{align}\label{decomp-tildezeta}
\tilde{\bzeta} = \sum\limits_{i=1}^N \sum\limits_{k\in \mathcal{M}_i}
\tilde{\bzeta}_{k,\Delta}^{(i)} + \sum\limits_{k=1}^M
\tilde{\bzeta}_{k,\Pi},~\tilde{\bzeta} = \tilde{\bu}, \tilde{\bv},
\end{align}
where $\tilde{\bzeta}_{k,\Delta}^{(i)} \in \bW_{k, \Delta}^{(i)}$, $\tilde{\bzeta}_{k,\Pi} \in \bW_{k,\Pi}$.

Using \eqref{def-T-k-i} and \eqref{decomp-tildezeta}, define a bilinear form on $\tilde{\bW}$ via
\begin{align}\label{def-operator-tilde-A}
\tilde{A}(\tilde{\bu}, \tilde{\bv}) = \sum\limits_{i=1}^N A_i(\tilde{\bu}^{(i)}, \tilde{\bv}^{(i)}),~~\forall
\tilde{\bu}, \tilde{\bv} \in \tilde{\bW},
\end{align}
where
\begin{align}\label{def-zeta}
\tilde{\bzeta}^{(i)} = \sum\limits_{k \in \mathcal{M}_i} (\tilde{\bzeta}_{k,\Delta}^{(i)} + T_{k}^{(i)}\tilde{\bzeta}_{k,\Pi}),~\tilde{\bzeta} = \tilde{\bu}, \tilde{\bv}.
\end{align}


Then, an SPD operator $\tilde{S}: \tilde{\bW}
\rightarrow \tilde{\bW}$ can be defined as follows
\begin{align}\label{def-operator-tilde-S}
(\tilde{S} \tilde{\bu}, \tilde{\bv}) = \tilde{A}(\tilde{\bu}, \tilde{\bv}),~~\forall \tilde{\bu},
\tilde{\bv} \in \tilde{\bW}.
\end{align}


Let $I_{\Gamma}$ be the natural injection from $\hat{\bW}$ into
$\tilde{\bW}$ (see \cite{BS2007}) such that for each $\bw_{\Pi} \in
\bW_{\Pi}$, $\bw_{k,\Delta} \in \bW_{k,\Delta}$ $(k=1,\cdots,M)$, we have
\begin{align}\label{def-natural-injection}
I_{\Gamma} \bw_{\Pi} = \bw_{\Pi},~ I_{\Gamma} \bw_{k,\Delta} =
T_{k}^{(i)} \bw_{k,\Delta} + T_{k}^{(j)}
\bw_{k,\Delta},
\end{align}
where $T_{k}^{(\nu)} (\nu = i,j)$ are defined in \eqref{def-T-k-i},
$i$ and
$j$ are the indices of the subdomains which satisfy $F_k = \partial
\Omega_i \cap
\partial \Omega_j$.

For any given $\hat{\bu}, \hat{\bv} \in \hat{\bW}$, using the decomposition \eqref{def-hatW} of $\hat{\bW}$, we have
\begin{align}\label{def-hatzeta}
\hat{\bzeta} = \sum\limits_{k=1}^M (\hat{\bzeta}_{k,\Delta} + \hat{\bzeta}_{k,\Pi}),~\hat{\bzeta} = \hat{\bu}, \hat{\bv} \in \hat{\bW},
\end{align}
where $\hat{\bzeta}_{k,\Delta} \in \bW_{k,\Delta}$, $\hat{\bzeta}_{k,\Pi} \in \bW_{k,\Pi}$.
From \eqref{def-natural-injection} and \eqref{def-hatzeta}, we find
\begin{align}\label{def-natural-injection-1-00}
 I_{\Gamma} \hat{\bzeta} = \sum\limits_{i=1}^N \sum\limits_{k\in \mathcal{M}_i} T_{k}^{(i)} \hat{\bzeta}_{k,\Delta}
 + \sum\limits_{k=1}^M \hat{\bzeta}_{k,\Pi},~ \hat{\bzeta} = \hat{\bu}, \hat{\bv}.
\end{align}
Furthermore, combining \eqref{def-AUV}, \eqref{def-operator-hat-S}
and \eqref{def-hatzeta}, yields
\begin{align}\label{def-AUV-eqiv}
(\hat{S}\hat{\bu}, \hat{\bv}) = A(\hat{\bu}, \hat{\bv}) = \sum\limits_{i=1}^N A_i(\hat{\bu}^{(i)},
\hat{\bv}^{(i)}),
\end{align}
where
\begin{align}\label{def-AUV-eqiv-i}
\hat{\bzeta}^{(i)} = \sum\limits_{k \in \mathcal{M}_i} T_{k}^{(i)}  (\hat{\bzeta}_{k,\Delta}  + \hat{\bzeta}_{k,\Pi}),~\hat{\bzeta} = \hat{\bu}, \hat{\bv}.
\end{align}
Using \eqref{def-operator-tilde-S}, \eqref{def-natural-injection-1-00} and \eqref{def-AUV-eqiv}, it is easy to derive that
\begin{align}\label{operator-hatS-tildeS-relation}
 \hat{S} = I_{\Gamma}^T \tilde{S} I_{\Gamma}.
\end{align}



For each interface $F_k = \partial \Omega_i \cap \partial \Omega_j$
($k=1,\cdots,M$), let the function vectors
\begin{align*}
\hat{\Phi}^{k} =
(\hat{\bphi}_1^k,\cdots,\hat{\bphi}_{n_k}^k)^T~\mbox{and}~\tilde{\Phi}^{k,\nu}
= (\tilde{\bphi}_1^{k,\nu},\cdots,\tilde{\bphi}_{n_k}^{k,\nu})^T,~\nu=i,j,
\end{align*}
and  $T^{(i,j)}: \bW_k^{(i)} \rightarrow \bW_k^{(j)}$ be the linear  basis transformation operator such that
\begin{align}\label{def-T-i-j}
T^{(i,j)} \tilde{\bphi}^{k,i}_l  = \tilde{\bphi}^{k,j}_l,~l=1,\cdots,n_k.
\end{align}

 Define the scaling operators $\check{D}_{F_{k}}^{(\nu)}$ $(\nu = i,j): \bU \rightarrow \bU$
$(\bU = \bW_{k}, \bW_{k}^{(i)}, \bW_{k}^{(j)})$ such that for any $\bw = \vec{w}^T \Psi$ with $\vec{w}
\in \mathbb{R}^{n_k}$ and $\Psi =
\hat{\Phi}^k,\tilde{\Phi}^{k,i}~\mbox{or}~ \tilde{\Phi}^{k,j}$, we have
\begin{align}\label{def-D-operator-matrix-check}
\check{D}_{F_{k}}^{(\nu)} \bw = \vec{w}^T (\vec{\check{D}}_{F_{k}}^{(\nu)})^T \Psi,
\end{align}
where $\vec{\check{D}}^{(\nu)}_{F_k}$ is a $n_k \times n_k$ scaling matrix, which satisfy
\begin{align}\label{def-D-trans-check}
\check{D}_{F_{k}}^{(i)} + \check{D}_{F_{k}}^{(j)}  = I,~\mbox{where $I$ is the identity operator,}
\end{align}
and
\begin{align}\label{def-D-trans-check-0}\nonumber
& A_i(\check{D}_{F_{k}}^{(j)} \bw,\check{D}_{F_{k}}^{(j)}\bw)
+ A_{j}(\check{D}_{F_{k}}^{(i)}T^{(i,j)} \bw, \check{D}_{F_{k}}^{(i)} T^{(i,j)} \bw) \\
& \le  A_i(D_{F_{k}}^{(j)} \bw,D_{F_{k}}^{(j)}\bw) +  A_{j}(D_{F_{k}}^{(i)} T^{(i,j)} \bw,
D_{F_{k}}^{(i)} T^{(i,j)} \bw),~~\forall \bw \in \bW_{k,\Delta}^{(i)},
\end{align}
here $T^{(i,j)}$ and  $D_{F_{k}}^{(\nu)}(\nu=i,j)$ are separately defined in \eqref{def-T-i-j} and \eqref{def-D-operator-matrix}.

The scaling matrices $\vec{\check{D}}_{F_{k}}^{(\nu)}$ $(\nu = i,j)$ in
\eqref{def-D-operator-matrix-check} are usually expressed as (see \cite{KCW2016})
\begin{align}\label{def-D-trans-check-2}
\vec{\check{D}}_{F_{k}}^{(\nu)} = \vec{T}_{F_k}^{-1} \vec{D}_{F_{k}}^{(\nu)} \vec{T}_{F_k},
\end{align}
where $\vec{D}_{F_{k}}^{(\nu)}$ and $\vec{T}_{F_k}$ are the matrices in \eqref{def-D-operator-matrix} and \eqref{def-bar-Phi-k-Delta1-nu}, respectively.


From \eqref{def-natural-injection} and \eqref{def-D-operator-matrix-check}, we can easily prove that the natural injection operator $I_{\Gamma}$ and the scaling operator $\check{D}_{F_k}^{(\nu)}$ $(\nu=i,j)$ satisfy the exchangeable property, namely
\begin{align}\label{property-I-D}
I_{\Gamma}  \check{D}_{F_k}^{(\nu)} \bw =  \check{D}_{F_k}^{(\nu)} I_{\Gamma}  \bw,~~\forall \bw \in \bW_{k,\Delta},~~\nu=i,j.
\end{align}


For each subdomain $\Omega_i (i=1,\cdots,N)$, let $T_{(i)}^k: \bW_k^{(i)} \rightarrow \bW_k (k\in \mathcal{M}_i)$ be the linear basis transformation operator such that
\begin{align}\label{def-T-i-k}
T_{(i)}^k \tilde{\bphi}^{k,i}_l = \hat{\bphi}^k_l,~l=1,\cdots,n_k.
\end{align}

By using $T_{(i)}^k$ and $\check{D}_{F_k}^{(i)} (k \in \mathcal{M}_i)$,
a linear operator $R_{\Delta,\Gamma}^{(i)}: \bW_{\Delta}^{(i)} \rightarrow \hat{\bW}$ can be defined by
\begin{align}\label{def-R-hatW-WDeltai-trans}
R_{\Delta,\Gamma}^{(i)} = \sum\limits_{k \in \mathcal{M}_i} \check{D}_{F_k}^{(i)} T_{(i)}^k R_{k,\Delta}^{(i)},
\end{align}
where $R_{k,\Delta}^{(i)}:\bW_{\Delta}^{(i)} \rightarrow \bW_{k,\Delta}^{(i)}$ is the restriction operator defined in  \eqref{restriction-operator}.


Using \eqref{def-R-hatW-WDeltai-trans}, we can easily verify that
\begin{align}\label{def-R-property-1}
R_{\Delta,\Gamma}^{(i)} \bw = \check{D}_{F_k}^{(i)} T_{(i)}^k \bw,~\forall \bw  \in \bW_{k,\Delta}^{(i)}, k \in \mathcal{M}_i, i=1,\cdots,N,
\end{align}
and for any given $\bw \in \bW_{k,\Delta}$ $(k=1,\cdots,M)$, from \eqref{def-R-property-1}, \eqref{def-T-k-i} and \eqref{def-D-trans-check}, we have
\begin{align}\label{def-R-hatW-W-Delta-i-new}
\bw = R_{\Delta,\Gamma}^{(i)} T_{k}^{(i)} \bw + R_{\Delta,\Gamma}^{(j)} T_{k}^{(j)} \bw,
\end{align}
where $i$ and $j$ are the indices of the subdomains which satisfy $F_k = \partial \Omega_i \cap \partial \Omega_j$.


The operator $R_{\Pi,\Gamma}:\bW_{\Pi} \rightarrow \tilde{\bW}$ is defined by
\begin{align}\label{R-WPi-barWPi-00}
R_{\Pi,\Gamma} = E_{\Pi} - \tilde{S}_{\Delta}^{-1} \tilde{S} E_{\Pi},
\end{align}
where
\begin{align}\label{def-tildeS-Delta-inv}
\tilde{S}_{\Delta}^{-1} = \sum\limits_{i=1}^N E_{\Delta}^{(i)} ((E_{\Delta}^{(i)})^T \tilde{S} E_{\Delta}^{(i)})^{-1} (E_{\Delta}^{(i)})^T,
\end{align}
here $E_{\Pi}:~\bW_{\Pi} \rightarrow \tilde{\bW}$ and $E_{\Delta}^{(i)}:~\bW_{\Delta}^{(i)} \rightarrow \tilde{\bW} (i=1,\cdots,N)$ are the interpolation operators defined in \eqref{interpolation-operator}. 


Using the above-mentioned preparations, we can present the adaptive BDDC preconditioned operator $ M_{BDDC}^{-1}: \hat{\bW} \rightarrow \hat{\bW}$ for solving the Schur
complement system \eqref{def-schur-system} in the variational form
as follows.


\begin{algorithm}(adaptive BDDC preconditioner)\label{algorithm-variation-BDDC-domain}~
Given ${\bg}\in \hat{\bW}$, the action $\bu_g = M_{BDDC}^{-1} {\bg}
\in \hat{\bW}$ is defined via the following four steps.
\begin{description}
\item[Step 1.] Find
\begin{align}\label{step2-ug}
\bu_{\Delta,a} = \sum\limits_{i=1}^N R_{\Delta,\Gamma}^{(i)} \bu^{\Delta,i}_a \in \hat{\bW},
\end{align}
where $\bu^{\Delta,i}_a \in \bW_{\Delta}^{(i)}(i=1,\cdots,N)$ such that
\begin{align}\label{step1-variation-form}
A_i(\bu^{\Delta,i}_a, \bv) = ( (R_{\Delta,\Gamma}^{(i)})^T {\bg}, \bv),~~\forall \bv \in \bW_{\Delta}^{(i)},
\end{align}
here the operator $R_{\Delta,\Gamma}^{(i)}$ is defined in \eqref{def-R-hatW-WDeltai-trans}.

\item[Step 2.]  Find $\bu_{\Pi} \in \bW_{\Pi}$ such that
\begin{align}\label{step3-variation-form-0001}
\tilde{A}(R_{\Pi,\Gamma} \bu_{\Pi}, R_{\Pi,\Gamma}\bv) = (\bg, \bv) - \tilde{A}(\sum\limits_{i=1}^{N} \bu^{\Delta,i}_a,\bv),~~\forall \bv\in
\bW_{\Pi},
\end{align}
where the operator $R_{\Pi,\Gamma}$ is defined in \eqref{R-WPi-barWPi-00}.

\item[Step 3.] Find
\begin{align}\label{step3-ug}
\bu_{\Delta,b} = \sum\limits_{i=1}^N R_{\Delta,\Gamma}^{(i)} \bu^{\Delta,i}_b \in \hat{\bW},
\end{align}
where  $\bu^{\Delta,i}_{b} \in \bW_{\Delta}^{(i)}(i=1,\cdots,N)$ satisfy that
\begin{align}\label{step3-variation-form-00}
A_i(\bu^{\Delta,i}_{b}, \bv) = - A_i(\bu_{\Pi}, \bv),~~\forall \bv \in \bW_{\Delta}^{(i)}.
\end{align}

\item[Step 4.] Compute
\begin{align}\label{step6-variation-form}
\bu_g = \bu_{\Delta,a} +  \bu_{\Pi} + \bu_{\Delta,b}.
\end{align}

\end{description}

\end{algorithm}

In the following,  we will derive the expressions for
$\bu_{\Delta,a}$, $\bu_{\Delta,b}$ and $\bu_{\Pi}$ in Algorithm
\ref{algorithm-variation-BDDC-domain}.


By using the definition of $\tilde{S}$ and \eqref{step1-variation-form}, we have
\begin{align}\label{equ-lemma-algorithm-M-BDDC-step1-proof-2}
\bu_a^{\Delta,i} = ((E_{\Delta}^{(i)})^T \tilde{S} E_{\Delta}^{(i)})^{-1} (R_{\Delta,\Gamma}^{(i)})^T {\bg}.
\end{align}

Combining \eqref{equ-lemma-algorithm-M-BDDC-step1-proof-2} and \eqref{step2-ug}, we obtain
\begin{align}\label{equ-lemma-algorithm-M-BDDC-step1}
\bu_{\Delta,a} = \sum\limits_{i=1}^N R_{\Delta,\Gamma}^{(i)} ((E_{\Delta}^{(i)})^T \tilde{S} E_{\Delta}^{(i)})^{-1} (R_{\Delta,\Gamma}^{(i)})^T {\bg}.
\end{align}


Similar to the derivation of \eqref{equ-lemma-algorithm-M-BDDC-step1}, we find
\begin{align}\label{equ-lemma-algorithm-M-BDDC-step3}
\bu_{\Delta,b} = -\sum\limits_{i=1}^{N} R_{\Delta,\Gamma}^{(i)} ((E_{\Delta}^{(i)})^T \tilde{S} E_{\Delta}^{(i)})^{-1}
(E_{\Delta}^{(i)})^T \tilde{S} E_{\Pi} \bu_{\Pi}.
\end{align}


Using \eqref{equ-lemma-algorithm-M-BDDC-step1-proof-2},  we have
\begin{align}\label{equ-lemma-algorithm-M-BDDC-step2-proof-2-001}
({\bg}, \bv) - \tilde{A}(\sum\limits_{i=1}^{N}\bu_a^{\Delta,i},\bv)
= ((\hat{E}_{\Pi})^T {\bg}, \bv) - (\sum\limits_{i=1}^{N} (E_{\Pi})^T \tilde{S} E_{\Delta}^{(i)}\bu^{\Delta,i}_a, \bv)
= (R_0 {\bg}, \bv),~\forall \bv\in
\bW_{\Pi}
\end{align}
 where
\begin{align}\label{equ-2-lemma-algorithm-M-BDDC}
R_0 &= (\hat{E}_{\Pi})^T - \sum\limits_{i=1}^N (E_{\Pi})^T \tilde{S} E_{\Delta}^{(i)} ((E_{\Delta}^{(i)})^T \tilde{S} E_{\Delta}^{(i)})^{-1} (R_{\Delta,\Gamma}^{(i)})^T,
\end{align}
here  $\hat{E}_{\Pi}: \bW_{\Pi} \rightarrow \hat{\bW}$ is the interpolation operator defined in \eqref{interpolation-operator}.

It is easy to know that 
\begin{align}\label{equ-lemma-algorithm-M-BDDC-step2-proof-1-001}
\tilde{A}(R_{\Pi,\Gamma}\bu_{\Pi},R_{\Pi,\Gamma}\bv)
=(F_{\Pi\Pi}\bu_{\Pi}, \bv)
,~\forall \bv\in \bW_{\Pi},
\end{align}
where
\begin{align}\label{def-FPiPi-2}
F_{\Pi\Pi} = (R_{\Pi,\Gamma})^T \tilde{S} R_{\Pi,\Gamma}.
\end{align}

Inserting  \eqref{equ-lemma-algorithm-M-BDDC-step2-proof-2-001} and \eqref{equ-lemma-algorithm-M-BDDC-step2-proof-1-001} into \eqref{step3-variation-form-0001}, and since $\bv\in \bW_{\Pi}$ is arbitrary, it implies
\begin{align}\label{equ-lemma-algorithm-M-BDDC-step2}
\bu_{\Pi} =  F_{\Pi\Pi}^{-1} R_0 {\bg}.
\end{align}

Substituting \eqref{equ-lemma-algorithm-M-BDDC-step1-proof-2}, \eqref{equ-lemma-algorithm-M-BDDC-step3} and \eqref{equ-lemma-algorithm-M-BDDC-step2} into \eqref{step6-variation-form},
we can arrive at the expression of $M_{BDDC}^{-1}$ defined in algorithm \ref{algorithm-variation-BDDC-domain} as follows:
\begin{align}\label{equ-1-lemma-algorithm-M-BDDC}
M_{BDDC}^{-1} = \sum\limits_{i=1}^N R_{\Delta,\Gamma}^{(i)} ((E_{\Delta}^{(i)})^T \tilde{S} E_{\Delta}^{(i)})^{-1} (R_{\Delta,\Gamma}^{(i)})^T + R_0^T F_{\Pi\Pi}^{-1} R_0.
\end{align}


In order to bound the condition number of $M_{BDDC}^{-1}$,
 we need to rewrite the preconditioned operator in a more concise form than \eqref{equ-1-lemma-algorithm-M-BDDC}.

Firstly, we give the equivalent expression of $F_{\Pi\Pi}$ defined
in \eqref{def-FPiPi-2}. For any given $\tilde{\bg}_{\Delta}\in
\tilde{\bW}_{\Delta}$, by using \eqref{def-tildeW-Delta}, there exists a decomposition
\begin{align}\label{def-tilde-g-00}
\tilde{\bg}_{\Delta} = \sum\limits_{i=1}^N E_{\Delta}^{(i)} \tilde{\bg}_{\Delta}^{(i)},~\mbox{where}~\tilde{\bg}_{\Delta}^{(i)} \in \bW_{\Delta}^{(i)}.
\end{align}
%
%

From \eqref{def-tildeS-Delta-inv} and \eqref{def-tilde-g-00}, one
has
\begin{align}\label{invtildeDeltaS-0}
\tilde{S}_{\Delta}^{-1} \tilde{S} \tilde{\bg}_{\Delta}
= \sum\limits_{i=1}^N \sum\limits_{j=1}^N E_{\Delta}^{(i)} ((E_{\Delta}^{(i)})^T \tilde{S} E_{\Delta}^{(i)})^{-1} (E_{\Delta}^{(i)})^T \tilde{S} E_{\Delta}^{(j)} \tilde{\bg}_{\Delta}^{(j)}.
\end{align}

Since the support property of the functions in
$\bW_{\Delta}^{(l)}(l=1,\cdots,N$), it implies that
\begin{align}\label{rel-deltai-deltaj}
(E_{\Delta}^{(i)})^T \tilde{S} E_{\Delta}^{(j)} = 0,~i\neq j.
\end{align}

%

Using \eqref{rel-deltai-deltaj} and \eqref{def-tilde-g-00}, we derive from \eqref{invtildeDeltaS-0} that
\begin{align}\label{def-FPiPi-2-proof-2}
\tilde{S}_{\Delta}^{-1} \tilde{S} \tilde{\bg}_{\Delta}
= \sum\limits_{i=1}^N E_{\Delta}^{(i)} ((E_{\Delta}^{(i)})^T \tilde{S} E_{\Delta}^{(i)})^{-1} (E_{\Delta}^{(i)})^T \tilde{S} E_{\Delta}^{(i)} \tilde{\bg}_{\Delta}^{(i)}
= \sum\limits_{i=1}^N E_{\Delta}^{(i)}\tilde{\bg}_{\Delta}^{(i)}
= \tilde{\bg}_{\Delta}.
\end{align}
Note that $\tilde{S}_{\Delta}^{-1} \tilde{S} \tilde{\bg} \in \tilde{\bW}_{\Delta}$, then \eqref{def-FPiPi-2-proof-2} implies
\begin{align}\label{def-FPiPi-2-proof-2-000}
(\tilde{S}_{\Delta}^{-1} \tilde{S})^2 = \tilde{S}_{\Delta}^{-1} \tilde{S}.
\end{align}

From \eqref{R-WPi-barWPi-00}, \eqref{def-FPiPi-2} and \eqref{def-FPiPi-2-proof-2-000},  together with the symmetry of $\tilde{S}$ and $\tilde{S}_{\Delta}^{-1}$, we get the equivalent form of $F_{\Pi\Pi}$ as follows
\begin{align}\label{def-FPiPi-2-proof-1}\nonumber
F_{\Pi\Pi} &= (E_{\Pi} - \tilde{S}_{\Delta}^{-1} \tilde{S} E_{\Pi})^T \tilde{S} (E_{\Pi} - \tilde{S}_{\Delta}^{-1} \tilde{S} E_{\Pi})\\\nonumber
&= (E_{\Pi})^T\tilde{S}(E_{\Pi} - \tilde{S}_{\Delta}^{-1} \tilde{S} E_{\Pi}) - (E_{\Pi})^T \tilde{S} (\tilde{S}_{\Delta}^{-1} \tilde{S}  - (\tilde{S}_{\Delta}^{-1} \tilde{S})^2) E_{\Pi}\\
&= (E_{\Pi})^T \tilde{S}E_{\Pi} - (E_{\Pi})^T\tilde{S} \tilde{S}_{\Delta}^{-1} \tilde{S} E_{\Pi}.
\end{align}

Next, we can derive the expression of the inverse operator of $\tilde{S}$ as
\begin{align}\label{def-inv-tildeS-2}\nonumber
\tilde{S}^{-1} &= \tilde{S}_{\Delta}^{-1}  + \tilde{S}_{\Delta}^{-1} \tilde{S} E_{\Pi} F_{\Pi\Pi}^{-1} (E_{\Pi})^T \tilde{S} \tilde{S}_{\Delta}^{-1}  - \tilde{S}_{\Delta}^{-1} \tilde{S} E_{\Pi} F_{\Pi\Pi}^{-1} (E_{\Pi})^T\\
&~~~ - E_{\Pi} F_{\Pi\Pi}^{-1} (E_{\Pi})^T \tilde{S}
\tilde{S}_{\Delta}^{-1}
  + E_{\Pi} F_{\Pi\Pi}^{-1} (E_{\Pi})^T.
\end{align}


In fact, according to the definition \eqref{def-tildeW} of
$\tilde{\bW}$, we only need to check that the operator
$B:=\tilde{S}^{-1}$ satisfies
\begin{align}\label{def-inv-tildeS-2-proof-1}
 B \tilde{S} \tilde{\bg}_{\Delta} = \tilde{\bg}_{\Delta},~\forall \tilde{\bg}_{\Delta}\in \tilde{\bW}_{\Delta},
\end{align}
and
\begin{align}\label{def-inv-tildeS-2-proof-2}
B \tilde{S} \tilde{\bg}_{\Pi} = \tilde{\bg}_{\Pi},~\forall \tilde{\bg}_{\Pi} \in \bW_{\Pi}.
\end{align}

From \eqref{def-FPiPi-2-proof-2} and \eqref{def-inv-tildeS-2}
, it is easy to know that \eqref{def-inv-tildeS-2-proof-1} is established.

%
By using \eqref{def-FPiPi-2-proof-1} and \eqref{def-inv-tildeS-2}, we find
\begin{align*}
B \tilde{S} \tilde{\bg}_{\Pi}
&= \tilde{S}_{\Delta}^{-1} \tilde{S} E_{\Pi} \tilde{\bg}_{\Pi}
+ \tilde{S}_{\Delta}^{-1} \tilde{S} E_{\Pi} F_{\Pi\Pi}^{-1} ((E_{\Pi})^T \tilde{S} \tilde{S}_{\Delta}^{-1} \tilde{S} E_{\Pi}
-  (E_{\Pi})^T \tilde{S} E_{\Pi}) \tilde{\bg}_{\Pi}\\\nonumber
&~~~ - E_{\Pi} F_{\Pi\Pi}^{-1} ( (E_{\Pi})^T \tilde{S} \tilde{S}_{\Delta}^{-1} \tilde{S} E_{\Pi}
 - (E_{\Pi})^T \tilde{S} E_{\Pi}) \tilde{\bg}_{\Pi}
 = \tilde{\bg}_{\Pi},~
\forall \tilde{\bg}_{\Pi} \in \bW_{\Pi},
\end{align*}
then \eqref{def-inv-tildeS-2-proof-2} holds.$\square$


In order to present  a concise form of $M_{BDDC}^{-1}$, we need
to introduce an averaging operator $E_D: \tilde{\bW} \rightarrow
\hat{\bW}$ as
\begin{align}\label{def-R-hatW-tildeW-form2}
E_D = \sum\limits_{i=1}^{N} R_{\Delta,\Gamma}^{(i)} R_{\Delta}^{(i)} +  R_{\Pi},
\end{align}
where the operator $R_{\Delta,\Gamma}^{(i)}$ is defined in \eqref{def-R-hatW-WDeltai-trans}, $R_{\Delta}^{(i)}:\tilde{\bW} \rightarrow \bW_{\Delta}^{(i)} (i=1,\cdots,N)$ and $R_{\Pi}:\tilde{\bW} \rightarrow \bW_{\Pi}$ are both the restriction operators defined in \eqref{restriction-operator}.


By using the definitions of the restriction operator and the interpolation operator, it is easy to verify that
\begin{align}\label{restriction-operator-property-000}
& R_{\Pi} E_{\Pi} = I,~~R_{\Delta}^{(i)} E_{\Pi} = 0,~~i=1,\cdots,N,\\\label{restriction-operator-property-001}
& R_{\Delta}^{(i)}  E_{\Delta}^{(j)} = \delta_{i,j}I,~~i,j=1,\cdots,N,
\end{align}
where $I$ is the identity operator.


Using the above-mentioned preparations, we have
\begin{theorem}
A concise form of the adaptive BDDC preconditioned operator
can be written as
\begin{align}\label{equ-1-lemma-algorithm-M-BDDC-form2}
M_{BDDC}^{-1} = E_D \tilde{S}^{-1} E_D^T,
\end{align}
where $\tilde{S}^{-1}$ and $E_D$ are defined in \eqref{def-inv-tildeS-2} and \eqref{def-R-hatW-tildeW-form2}, respectively.
\end{theorem}
\begin{proof}
From \eqref{def-R-hatW-tildeW-form2} and the interpolation operator $\hat{E}_{\Pi}$, we have
\begin{align}\label{M-BDDC-form-2-proof-0}\nonumber
E_D \tilde{S}^{-1} E_D^T 
&= \left(\sum\limits_{i=1}^{N} R_{\Delta,\Gamma}^{(i)} R_{\Delta}^{(i)} + \hat{E}_{\Pi} R_{\Pi} \right) \tilde{S}^{-1} \left(\sum\limits_{j=1}^{N} (R_{\Delta}^{(j)})^T (R_{\Delta,\Gamma}^{(j)})^T + (R_{\Pi})^T (\hat{E}_{\Pi})^T\right)\\
&= M_1 + M_2 + M_2^T + M_3,
\end{align}
where
\begin{align*}
& M_1 = \sum\limits_{i=1}^{N} R_{\Delta,\Gamma}^{(i)} R_{\Delta}^{(i)} \tilde{S}^{-1} \sum\limits_{j=1}^{N} (R_{\Delta}^{(j)})^T (R_{\Delta,\Gamma}^{(j)})^T, ~~M_2 = \sum\limits_{i=1}^{N} R_{\Delta,\Gamma}^{(i)} R_{\Delta}^{(i)} \tilde{S}^{-1} (R_{\Pi})^T (\hat{E}_{\Pi})^T,\\
& M_3 = \hat{E}_{\Pi} R_{\Pi} \tilde{S}^{-1} (R_{\Pi})^T (\hat{E}_{\Pi})^T.
\end{align*}

Using the above expression of $M_1$, together with \eqref{def-inv-tildeS-2} and the second equation of    \eqref{restriction-operator-property-000}, we obtain
\begin{align}\label{M-BDDC-form-2-proof-04-000}
M_1 = E \sum\limits_{j=1}^{N} (R_{\Delta}^{(j)})^T (R_{\Delta,\Gamma}^{(j)})^T + E \tilde{S} E_{\Pi} F_{\Pi\Pi}^{-1} (E_{\Pi})^T \tilde{S} E^T,
\end{align}
where
\begin{align}\label{expression-E}
E = \sum\limits_{i=1}^{N} R_{\Delta,\Gamma}^{(i)} R_{\Delta}^{(i)} \tilde{S}_{\Delta}^{-1}.
\end{align}


Inserting \eqref{def-tildeS-Delta-inv} into the \eqref{expression-E} and using \eqref{restriction-operator-property-001}, we can rewrite the above operator $E$ as
\begin{align}\label{M-BDDC-form-2-proof-01-00}
E = \sum\limits_{i=1}^{N} \sum\limits_{j=1}^N  R_{\Delta,\Gamma}^{(i)} R_{\Delta}^{(i)}  E_{\Delta}^{(j)} ((E_{\Delta}^{(j)})^T \tilde{S} E_{\Delta}^{(j)})^{-1} (E_{\Delta}^{(j)})^T
= \sum\limits_{i=1}^N R_{\Delta,\Gamma}^{(i)} ((E_{\Delta}^{(i)})^T \tilde{S} E_{\Delta}^{(i)})^{-1} (E_{\Delta}^{(i)})^T.
\end{align}


Similarly, we can get the expressions of $M_2$ and $M_3$ as follows 
\begin{align}\label{M-BDDC-form-2-proof-05}
M_2 = - E \tilde{S} E_{\Pi} F_{\Pi\Pi}^{-1} (\hat{E}_{\Pi})^T, ~M_3 = \hat{E}_{\Pi} F_{\Pi\Pi}^{-1} (\hat{E}_{\Pi})^T.
\end{align}

Further, substituting \eqref{M-BDDC-form-2-proof-04-000} and \eqref{M-BDDC-form-2-proof-05} into \eqref{M-BDDC-form-2-proof-0}, and using \eqref{M-BDDC-form-2-proof-01-00}, \eqref{restriction-operator-property-001} and \eqref{equ-2-lemma-algorithm-M-BDDC}, we have
\begin{align*}\nonumber
E_D \tilde{S}^{-1} E_D^T
&= E \sum\limits_{j=1}^{N} (R_{\Delta}^{(j)})^T (R_{\Delta,\Gamma}^{(j)})^T
   + \left(\hat{E}_{\Pi} - E \tilde{S} E_{\Pi} \right)  F_{\Pi\Pi}^{-1} \left((\hat{E}_{\Pi})^T - (E_{\Pi})^T \tilde{S} E^T\right)\\
&= \sum\limits_{i=1}^N R_{\Delta,\Gamma}^{(i)} ((E_{\Delta}^{(i)})^T \tilde{S} E_{\Delta}^{(i)})^{-1} (R_{\Delta,\Gamma}^{(i)})^T
   + R_0^T F_{\Pi\Pi}^{-1} R_0,
\end{align*}
this combines with \eqref{equ-1-lemma-algorithm-M-BDDC}, we complete the proof of \eqref{equ-1-lemma-algorithm-M-BDDC-form2}.

\end{proof}

Using \eqref{operator-hatS-tildeS-relation} and \eqref{equ-1-lemma-algorithm-M-BDDC-form2}, we present the preconditioned systems of \eqref{def-schur-system} as
\begin{align}\label{oper-hatG}
\hat{G} = E_D \tilde{S}^{-1} E_D^T I_{\Gamma}^T \tilde{S} I_{\Gamma}.
\end{align}
We will give bounds on the condition number of $\hat{G}$ in the next section.

\section{Analysis of  the condition number}\label{sec:5}
\setcounter{equation}{0}

First of all, we want to estimate the minimum eigenvalue of $\hat{G}$. For this reason, we firstly derive the following partition of unity condition
\begin{align}\label{lemma-R-hatW-tildeW-trans}
E_D I_{\Gamma}  = I,
\end{align}
where the natural injection operator $I_{\Gamma}$ and the averaging operator $E_D$ are separately
defined in \eqref{def-natural-injection} and \eqref{def-R-hatW-tildeW-form2},
$I$ is the identity operator.


In fact, for any $\hat{\bu} \in \hat{\bW}$, by using \eqref{def-hatW}, we have the decomposition
\begin{align}\label{decomposition-hatw}
\hat{{\bu}} = \sum\limits_{k=1}^M \hat{\bu}_{k,\Delta} + \hat{\bu}_{\Pi},
\end{align}
where $\hat{\bu}_{k,\Delta} \in \bW_{k,\Delta}$, $\hat{\bu}_{\Pi}\in \bW_{\Pi}$.


From the definition of  $I_{\Gamma}$ and \eqref{decomposition-hatw}, we have
\begin{align}\label{def-I-hatW-tildeW-hatw}
I_{\Gamma} \hat{\bu} = \sum\limits_{i=1}^N \sum\limits_{k \in \mathcal{M}_i} T_k^{(i)} \hat{\bu}_{k,\Delta} + \hat{\bu}_{\Pi},
\end{align}
where $T_k^{(i)}$ is defined in \eqref{def-T-k-i}.


By the definitions of $E_D$, $R_{\Delta}^{(i)}(i=1,\cdots,N)$ and $R_{\Pi}$,
together with \eqref{def-I-hatW-tildeW-hatw}, \eqref{decomposition-hatw} and
the property \eqref{def-R-hatW-W-Delta-i-new} of $R_{\Delta,\Gamma}^{(l)}(l=1,\cdots,N)$, it follows that
\begin{align*} %
E_D I_{\Gamma} \hat{{\bu}}
&= \sum\limits_{i=1}^{N} R_{\Delta,\Gamma}^{(i)} \sum\limits_{k \in \mathcal{M}_i} T_k^{(i)} \hat{\bu}_{k,\Delta} + \hat{\bu}_{\Pi}\\\nonumber
&= \sum\limits_{k=1}^{M}(R_{\Delta,\Gamma}^{(i)} T_{k}^{(i)} \hat{\bu}_{k,\Delta} + R_{\Delta,\Gamma}^{(j)} T_{k}^{(j)} \hat{\bu}_{k,\Delta}) + \hat{\bu}_{\Pi}\\\nonumber
&= \sum\limits_{k=1}^{M} \hat{\bu}_{k,\Delta} + \hat{\bu}_{\Pi}
= \hat{{\bu}}
\end{align*}
where we have used the assumption that each interface $F_k = \partial \Omega_i \cap \partial \Omega_j$ in the second equality.

From this and note that $\hat{\bu} \in \hat{\bW}$ is arbitrary, the proof of \eqref{lemma-R-hatW-tildeW-trans} is completed.$\square$


By using \eqref{lemma-R-hatW-tildeW-trans}, an argument similar to Lemma 3.4 in \cite{BS2007}
shows that the minimum eigenvalue of the preconditioned system  $\hat{G}$ satisfies
\begin{align}\label{lambdamin-hatG-oper}
\lambda_{\min}(\hat{G}) \ge 1.
\end{align}


Then, we are in the position to derive an upper bound for the maximum eigenvalue of $\hat{G}$.
Let $P_D: \tilde{\bW} \rightarrow \tilde{\bW}$ be the jump operator defined by
\begin{align}\label{def-Pd-0}
P_D = I - I_{\Gamma} E_D.
\end{align}


A conversion process similar to the maximum eigenvalue of $\hat{G}$ in an algebraic framework (see \cite{KC2015}, Lemma 3.1 and Lemma 3.2) shows that
\begin{align}\label{transform-1}
\lambda_{\max}(\hat{G}) \le \lambda_{\max}(G_d),
\end{align}
where the operator $G_d = P_D^T \tilde{S} P_D \tilde{S}^{-1}$.


Note that $G_d$ and $\tilde{S}^{-1} P_D^T \tilde{S} P_D$ share the same set of nonzero eigenvalues.
From this and using \eqref{transform-1},  the symmetry of $\tilde{S}^{-1} P_D^T \tilde{S} P_D$ with respect to the bilinear form $\tilde{A}(\cdot,\cdot)$ and the definition \eqref{def-operator-tilde-S} of $\tilde{S}$, we can arrive at
\begin{align}\label{5-26-proof-2-000}
\lambda_{\max}(\hat{G}) \le \max\limits_{\tilde{\bw} \in  \tilde{\bW} \backslash \{{\bf 0}\}} \frac{\tilde{A}(\tilde{S}^{-1} P_D^T \tilde{S} P_D \tilde{\bw}, \tilde{\bw})}{\tilde{A}(\tilde{\bw}, \tilde{\bw})} = \max\limits_{\tilde{\bw} \in  \tilde{\bW} \backslash \{{\bf 0}\}} \frac{\tilde{A}(P_D \tilde{\bw}, P_D \tilde{\bw})}{\tilde{A}(\tilde{\bw}, \tilde{\bw})}.
\end{align}


For any given $\tilde{\bw} \in \tilde{\bW}$, we can derive the decomposition formula of $P_D \tilde{\bw}$.
Using the definition \eqref{def-tildeW} of $\tilde{\bW}$, we have
\begin{align}\label{def-PD-tildew}
\tilde{\bw} = \sum\limits_{i=1}^N \sum\limits_{k \in \mathcal{M}_i} \bw_{k,\Delta}^{(i)} + \bw_{\Pi},~\mbox{where}~\bw_{k,\Delta}^{(i)} \in \bW_{k,\Delta}^{(i)}, \bw_{\Pi} = \sum\limits_{k=1}^M \bw_{k,\Pi} \in \bW_{\Pi}.
\end{align}


By the definitions of $P_D$, $E_D$, $R_{\Pi}$ and $I_{\Gamma}$, and using the decomposition \eqref{def-PD-tildew},
$P_D \tilde{\bw}$ can be rewrite as follows:
\begin{align}\label{lemmaLT-I-equ-proof-1-001-0}\nonumber
P_D \tilde{\bw} &= \tilde{\bw} - I_{\Gamma}E_D \tilde{\bw}\\\nonumber
&=\tilde{\bw} - I_{\Gamma}(\sum\limits_{i=1}^{N} R_{\Delta,\Gamma}^{(i)} R_{\Delta}^{(i)} +  R_{\Pi}) \tilde{\bw}\\\nonumber
&= \sum\limits_{i=1}^N \sum\limits_{k \in \mathcal{M}_i} \bw_{k,\Delta}^{(i)} + \bw_{\Pi} - I_{\Gamma}\sum\limits_{i=1}^{N} R_{\Delta,\Gamma}^{(i)} R_{\Delta}^{(i)} \tilde{\bw} -\bw_{\Pi} \\
&= \sum\limits_{i=1}^N \sum\limits_{k \in \mathcal{M}_i} \bw_{k,\Delta}^{(i)}
   - I_{\Gamma}\sum\limits_{i=1}^{N} R_{\Delta,\Gamma}^{(i)} R_{\Delta}^{(i)} \tilde{\bw}
\end{align}

By using \eqref{def-PD-tildew}, the definitions of $R_{\Delta}^{(i)} (i=1,\cdots,N)$ and $I_{\Gamma}$, the properties \eqref{def-R-property-1} and  \eqref{property-I-D},
we find the second term in \eqref{lemmaLT-I-equ-proof-1-001-0} satisfies
\begin{align}\label{equ-IR-1}\nonumber
I_{\Gamma} \sum\limits_{i=1}^{N} R_{\Delta,\Gamma}^{(i)} R_{\Delta}^{(i)} \tilde{\bw}
&= I_{\Gamma} \sum\limits_{i=1}^{N} R_{\Delta,\Gamma}^{(i)} \sum\limits_{k \in \mathcal{M}_i} \bw_{k,\Delta}^{(i)} \\\nonumber
&= \sum\limits_{i=1}^{N} \sum\limits_{k \in \mathcal{M}_i} \check{D}_{F_k}^{(i)} I_{\Gamma} T_{(i)}^k \bw_{k,\Delta}^{(i)} \\
&= \sum\limits_{i=1}^{N} \sum\limits_{k \in \mathcal{M}_i} \check{D}_{F_k}^{(i)} (\bw_{k,\Delta}^{(i)} + T^{(i,j)} \bw_{k,\Delta}^{(i)})
\end{align}
where we assume that $F_k = \partial \Omega_i \cap \partial \Omega_j (k \in \mathcal{M}_i, 1\le i\le N)$,
and the basis transformation operators $T_{(i)}^k$ and $T^{(i,j)}$ are separately defined in \eqref{def-T-i-k} and \eqref{def-T-i-j}.

Substituting \eqref{equ-IR-1} into \eqref{lemmaLT-I-equ-proof-1-001-0}, and using the property of $\check{D}_{F_k}^{(\nu)}(\nu=i,j)$ in \eqref{def-D-trans-check}, we can obtain the decomposition of $P_D \tilde{\bw}$ as follows
\begin{align}\label{lemmaLT-I-equ-proof-1-001}\nonumber
P_D \tilde{\bw} &= \sum\limits_{i=1}^N \sum\limits_{k \in \mathcal{M}_i} \left((\check{D}_{F_k}^{(i)} + \check{D}_{F_k}^{(j)})\bw_{k,\Delta}^{(i)} - \check{D}_{F_k}^{(i)} (\bw_{k,\Delta}^{(i)} + T^{(i,j)} \bw_{k,\Delta}^{(i)}) \right)\\\nonumber
&= \sum\limits_{i=1}^N \sum\limits_{k \in \mathcal{M}_i} (\check{D}_{F_k}^{(j)} \bw_{k,\Delta}^{(i)} - \check{D}_{F_k}^{(i)} T^{(i,j)} \bw_{k,\Delta}^{(i)})\\
&= \sum\limits_{i=1}^N \sum\limits_{k \in \mathcal{M}_i}\check{D}_{F_k}^{(j)}(\bw_{k,\Delta}^{(i)} - T^{(j,i)} \bw_{k,\Delta}^{(j)})
\end{align}


Using the decomposition \eqref{lemmaLT-I-equ-proof-1-001}, we can derive the following lemma:
\begin{lemma}\label{lemma-max}
For a given tolerance $\Theta \ge 1$, the maximum eigenvalue of the adaptive BDDC preconditioned system $\hat{G}$ satisfis
\begin{align}\label{eig-max}
\lambda_{\max}(\hat{G}) \le C \Theta,
\end{align}
where $C = 2 C_F^2$ and $C_F = \max\limits_{i}\{f_i\}$, here $f_i$ denotes the number of interface on $\partial \Omega_i$. 
\end{lemma}


According to \eqref{5-26-proof-2-000}, in order to give the proof of \eqref{eig-max},  we only need to show that
\begin{align*}
\max\limits_{\tilde{\bw} \in  \tilde{\bW} \backslash \{{\bf 0}\}} \frac{\tilde{A}(P_D \tilde{\bw}, P_D \tilde{\bw})}{\tilde{A}(\tilde{\bw}, \tilde{\bw})} \le C \Theta.
\end{align*}
In view of the definitions \eqref{def-operator-tilde-A}, \eqref{def-zeta} of the bilinear form $\tilde{A}(\cdot,\cdot)$, and the decompositions \eqref{def-PD-tildew} and \eqref{lemmaLT-I-equ-proof-1-001}, it is equivalent to show that
\begin{align}\label{Pb-max-eigenvalue-estimate-equiv-1-000-a}
\sum\limits_{i=1}^N A_i((P_D \tilde{\bw})^{(i)}, (P_D \tilde{\bw})^{(i)}) \le C
\Theta \sum\limits_{i=1}^N A_i(\tilde{\bw}^{(i)},\tilde{\bw}^{(i)}),
\end{align}
where
\begin{align}\label{PDi}
(P_D\tilde{\bw})^{(i)} = \sum\limits_{k \in \mathcal{M}_i} \check{D}_{F_k}^{(j)}(\bw_{k,\Delta}^{(i)} - \tilde{\bw}_{k,\Delta}^{(i)}),~
\tilde{\bw}^{(i)} = \sum\limits_{k \in \mathcal{M}_i} (\bw_{k,\Delta}^{(i)} + \bw_{k,\Pi}^{(i)}),
\end{align}
here
$$\tilde{\bw}_{k,\Delta}^{(i)} = T^{(j,i)} \bw_{k,\Delta}^{(j)},~~\bw_{k,\Pi}^{(i)} = T_k^{(i)} \bw_{k,\Pi}.$$


Firstly, by \eqref{PDi} and the essential properties \eqref{def-D-trans-check-0}, \eqref{def-TFkDelta-2-equ-1}, we have 
\begin{align}\label{final-ineq-1}\nonumber
&\sum\limits_{i=1}^N A_i((P_D\tilde{\bw})^{(i)},(P_D\tilde{\bw})^{(i)})\\\nonumber
&=\sum\limits_{i=1}^N A_i(\sum\limits_{k \in \mathcal{M}_i}\check{D}_{F_k}^{(j)}(\bw_{k,\Delta}^{(i)} - \tilde{\bw}_{k,\Delta}^{(i)}),\sum\limits_{k \in \mathcal{M}_i}\check{D}_{F_k}^{(j)}(\bw_{k,\Delta}^{(i)} - \tilde{\bw}_{k,\Delta}^{(i)}))\\\nonumber
&\le 2 C_F \sum\limits_{i=1}^N\sum\limits_{k \in \mathcal{M}_i}\left(A_i(\check{D}_{F_k}^{(j)} \bw_{k,\Delta}^{(i)}, \check{D}_{F_k}^{(j)} \bw_{k,\Delta}^{(i)})
+ A_i(\check{D}_{F_k}^{(j)} \tilde{\bw}_{k,\Delta}^{(i)}, \check{D}_{F_k}^{(j)} \tilde{\bw}_{k,\Delta}^{(i)})\right)\\\nonumber
&= 2 C_F \sum\limits_{i=1}^N \sum\limits_{k \in \mathcal{M}_i}
\left(A_i(\check{D}_{F_k}^{(j)} \bw_{k,\Delta}^{(i)},\check{D}_{F_k}^{(j)}\bw_{k,\Delta}^{(i)})
+ A_j(\check{D}_{F_k}^{(i)}\tilde{\bw}_{k,\Delta}^{(j)}, \check{D}_{F_k}^{(i)}\tilde{\bw}_{k,\Delta}^{(j)})\right)\\\nonumber
&\le  2 C_F \sum\limits_{i=1}^N \sum\limits_{k \in \mathcal{M}_i}
\left(A_i(D_{F_k}^{(j)} \bw_{k,\Delta}^{(i)},D_{F_k}^{(j)}\bw_{k,\Delta}^{(i)}) +  A_j(D_{F_k}^{(i)}\tilde{\bw}_{k,\Delta}^{(j)}, D_{F_k}^{(i)}\tilde{\bw}_{k,\Delta}^{(j)})\right)\\
&\le 2 C_F \Theta \sum\limits_{i=1}^N \sum\limits_{k \in \mathcal{M}_i}
A_i(\bar{\bw}_{k,\Delta}^{(i)} + \bar{\bw}_{k,\Pi}^{(i)}, \bar{\bw}_{k,\Delta}^{(i)} + \bar{\bw}_{k,\Pi}^{(i)})
\end{align}
where
\begin{align}\label{step-1-param-property}
\bar{\bw}_{k,\Delta}^{(i)} = \bar{T}_k^{(i)} \bw_{k,\Delta}^{(i)} \in \bar{\bW}_{k,\Delta}^{(i)},~~\bar{\bw}_{k,\Pi}^{(i)}  = \bar{T}_k^{(i)} \bw_{k,\Pi}^{(i)} \in \bar{\bW}_{k,\Pi}^{(i)},
\end{align}
here the linear basis transformation operator $\bar{T}^{(i)}_k:\bW_k^{(i)} \rightarrow \bar{\bW}_k^{(i)}$ is defined by
\begin{align*}
\bar{T}_k^{(i)} \tilde{\bphi}^{k,i}_l = \tilde{\bar{\bphi}}^{k,i}_l,~l=1,\cdots,n_k,
\end{align*}
and the basis functions $\{\tilde{\bar{\bphi}}^{k,i}_l\}_{l=1}^{n_k}$ is defined in \eqref{basis-2}.

Secondly, for each $k \in \mathcal{M}_i$, by using \eqref{PDi}, we obtain
\begin{align}\label{w-i-1}\nonumber
\tilde{\bw}^{(i)}
&= \bw_{k,\Delta}^{(i)} +  \bw_{k,\Pi}^{(i)} + \sum\limits_{m \in \mathcal{M}_i \atop m \neq k}(\bw_{m,\Delta}^{(i)} + \bw_{m,\Pi}^{(i)})\\\nonumber
&= (\bar{\bw}_{k,\Delta}^{(i)} - (\bar{\bw}_{k,\Delta}^{(i)} - \bw_{k,\Delta}^{(i)}))  + (\bar{\bw}_{k,\Pi}^{(i)} - (\bar{\bw}_{k,\Pi}^{(i)} - \bw_{k,\Pi}^{(i)})) + \sum\limits_{m \in \mathcal{M}_i \atop  m \neq k}(\bw_{m,\Delta}^{(i)} + \bw_{m,\Pi}^{(i)})\\
&= (\bar{\bw}_{k,\Delta}^{(i)} + \bar{\bw}_{k,\Pi}^{(i)}) + \bw_1 + \bw_2
\end{align}
where $\bw_1=  \sum\limits_{m \in \mathcal{M}_i \atop m \neq k} (\bw_{m,\Delta}^{(i)} + \bw_{m,\Pi}^{(i)})$ and $\bw_2= - (\bar{\bw}_{k,\Delta}^{(i)} - \bw_{k,\Delta}^{(i)}) - (\bar{\bw}_{k,\Pi}^{(i)} - \bw_{k,\Pi}^{(i)})$.


Obviously,
\begin{align}\label{subset-equ-1}
\bar{\bw}_{k,\Delta}^{(i)} + \bar{\bw}_{k,\Pi}^{(i)} \in \bar{\bW}_k^{(i)},~\bw_1 \in  \oplus_{m \in \mathcal{M}_i \atop m \neq k} \bW_m^{(i)},
\end{align}
and from \eqref{step-1-param-property} and the definition \eqref{def-Z-k-nu} of $\bZ_k^{(i)}$, we know that
\begin{align}\label{subset-equ-2}
\bw_2 \in  \bZ_k^{(i)}.
\end{align}
Using \eqref{w-i-1}, \eqref{subset-equ-1}, \eqref{subset-equ-2} and the orthogonality condition \eqref{def-w1-property-1-000}, we have
\begin{align}\label{4-1-2-tildew}\nonumber
A_i(\tilde{\bw}^{(i)},\tilde{\bw}^{(i)})
&= A_i(\bar{\bw}_{k,\Delta}^{(i)} + \bar{\bw}_{k,\Pi}^{(i)} + \bw_1 + \bw_2,\bar{\bw}_{k,\Delta}^{(i)} + \bar{\bw}_{k,\Pi}^{(i)} + \bw_1 + \bw_2) \\\nonumber
&=A_i(\bar{\bw}_{k,\Delta}^{(i)} + \bar{\bw}_{k,\Pi}^{(i)},\bar{\bw}_{k,\Delta}^{(i)} + \bar{\bw}_{k,\Pi}^{(i)}) + A_i(\bw_1 + \bw_2, \bw_1 + \bw_2)\\
&\ge A_i(\bar{\bw}_{k,\Delta}^{(i)} + \bar{\bw}_{k,\Pi}^{(i)},\bar{\bw}_{k,\Delta}^{(i)} + \bar{\bw}_{k,\Pi}^{(i)})
\end{align}
for all $k \in \mathcal{M}_i$.

Finally, the estimate \eqref{Pb-max-eigenvalue-estimate-equiv-1-000-a} follows from \eqref{final-ineq-1} and \eqref{4-1-2-tildew}.


By the results of \eqref{lambdamin-hatG-oper} and Lemma \ref{lemma-max}, we can obtain the following theorem.
\begin{theorem}\label{condition-number}
For a given tolerance $\Theta \ge 1$, we obtain the following condition number bound of the adaptive BDDC preconditioned systems $\hat{G}$ satisfying
\begin{align}\label{eig-condition}
\kappa(\hat{G}) \le C \Theta,
\end{align}
where  $C$ is a constant which is just depending on the maximum number of interfaces per each subdomain.
\end{theorem}


\section{Numerical results}\label{sec:6}\setcounter{equation}{0}

\begin{figure}[H]
  \centering
  \includegraphics[width=0.25\textwidth]{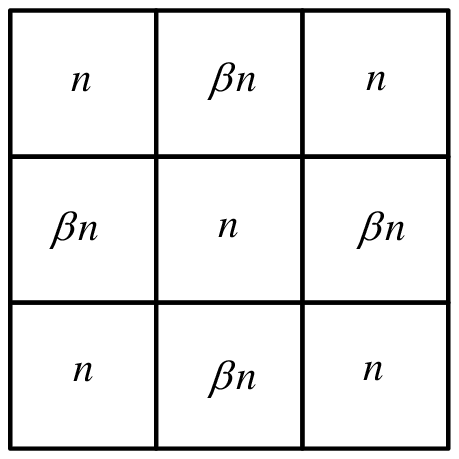}\\
  \caption{The checkerboard distribution of the local problem size in a geometrically conforming partitions.}\label{fig-checkerboard}
\end{figure}

In this section, we will present some numerical results of our adaptive BDDC algorithm for solving the Schur complement system \eqref{def-schur-system}.
We set the zero-order coefficient $\varepsilon = 1$ and the given region $\Omega = (0, 1)^2$ is decomposed into $N$ geometrically conforming or unconforming square subdomains.
Each subdomain is divided into a uniform triangulation mesh with $n$ or $\beta n$ elements in each direction distributed as checkerboard ($\beta \neq 1$ means the grids is non-matching),
the case with geometrically conforming subdomains see Figure \ref{fig-checkerboard}.
The PCG method is stopped when the relative residual  is reduced by the factor of $10^{-10}$.
For each interface, we emphasize that the nonmortar side is the one whose domain has larger step size.

In our adaptive BDDC algorithm, we set the tolerance $\Theta = 1+ \min\{log(n), log(\beta n)\}$ for a given mesh partition, the transformation matrix $\vec{T}_{F_k}$ and the scaling matricies $\vec{\check{D}}_{F_k}^{(\nu)} (\nu = i,j)$ in each interface $F_k$ are defined in \eqref{def-T-1} and \eqref{def-D-trans-check-2}, respectively.
Therefore, the algorithm is uniquely determined by another scaling matrices $\vec{D}_{F_k}^{(\nu)}(\nu=i,j)$ in each interface.
In the following experiments, we separately denote the algorithm with $\vec{D}_{F_k}^{(\nu)}(\nu=i,j, k=1,\cdots,M)$ defined in \eqref{DiFD} and \eqref{DiFD-deluxe} as M1 and M2.
we will investigate the robustness of these methods by some important parameters, such as Iter(number of iterations), $\lambda_{\min}$(minimum eigenvalue), $\lambda_{\max}$(maximum eigenvalue), $\kappa$(condition number), pnum(number of primal unknowns), ppnum(proportion of the total number of primal unknowns to the total number of dofs).

Firstly, we present some numerical results with geometrically conforming square subdomains.
Without loss of generality we assume that the mesh parameter $\beta \neq 1$ and  the space $X_h$ is associated with the $\mathcal{P}_2$ Lagrange finite element space.

\begin{example}\label{example-1}
Consider model problem \eqref{continuous-var} with $\rho(x) = 1$ for all $\Omega_i$.
\end{example}

In Table \ref{table-1-A}, the results are presented by increasing $n$ and with a fixed subdomain partition ($N = 3^2$) for Example \ref{example-1} with the mesh parameter $\beta = 1/2$.
We can observe that the total number of primal unknowns are the same in both methods and independent of $n$, and M2 has
lesser iterations than M1.

\begin{table}[H]
\centering\caption{Performance of the two methods for Example \ref{example-1}}
\label{table-1-A}\vskip 0.1cm
\begin{tabular}{{|c|c|c|c|c|c|}}\hline
n                  &method   &Iter  & $\lambda_{\min}$   & $\lambda_{\max}$  &pnum \\\hline
\multirow{2}{*}{12}  &M1  &9   &1.0014   &1.5148    &   16  \\
                     &M2  &6   &1.0001   &1.3076    &   16  \\\hline
\multirow{2}{*}{24}  &M1  &9   &1.0018   &1.6696    &   16  \\
                     &M2  &6   &1.0000   &1.4564    &   16  \\\hline
\multirow{2}{*}{48}  &M1  &9   &1.0024   &1.8275    &   16  \\
                     &M2  &7   &1.0000   &1.6177    &   16   \\\hline
\end{tabular}
\end{table}

\begin{figure}[H]
\begin{minipage}[t]{0.49\linewidth}
    \includegraphics[width=0.9\textwidth]{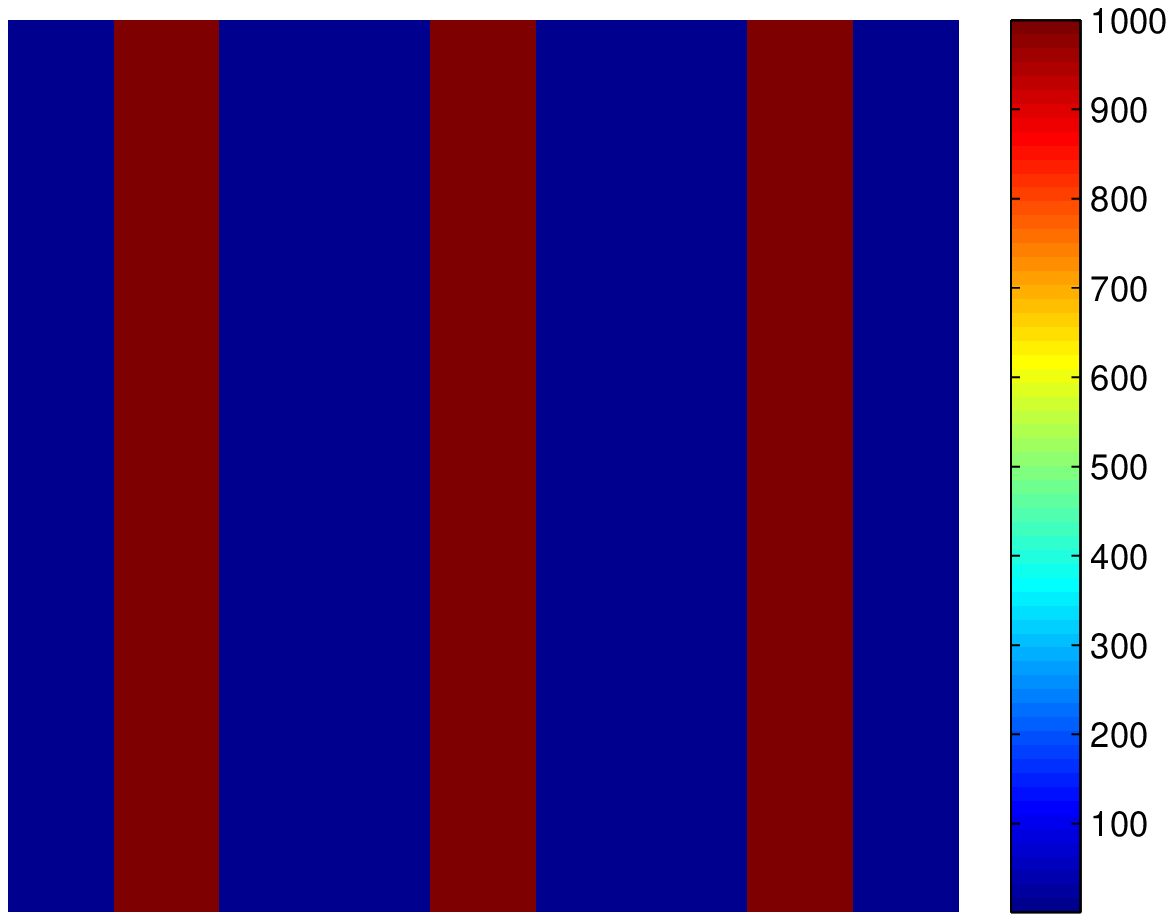}
\end{minipage}
\begin{minipage}[t]{0.49\linewidth}
    \includegraphics[width=0.9\textwidth]{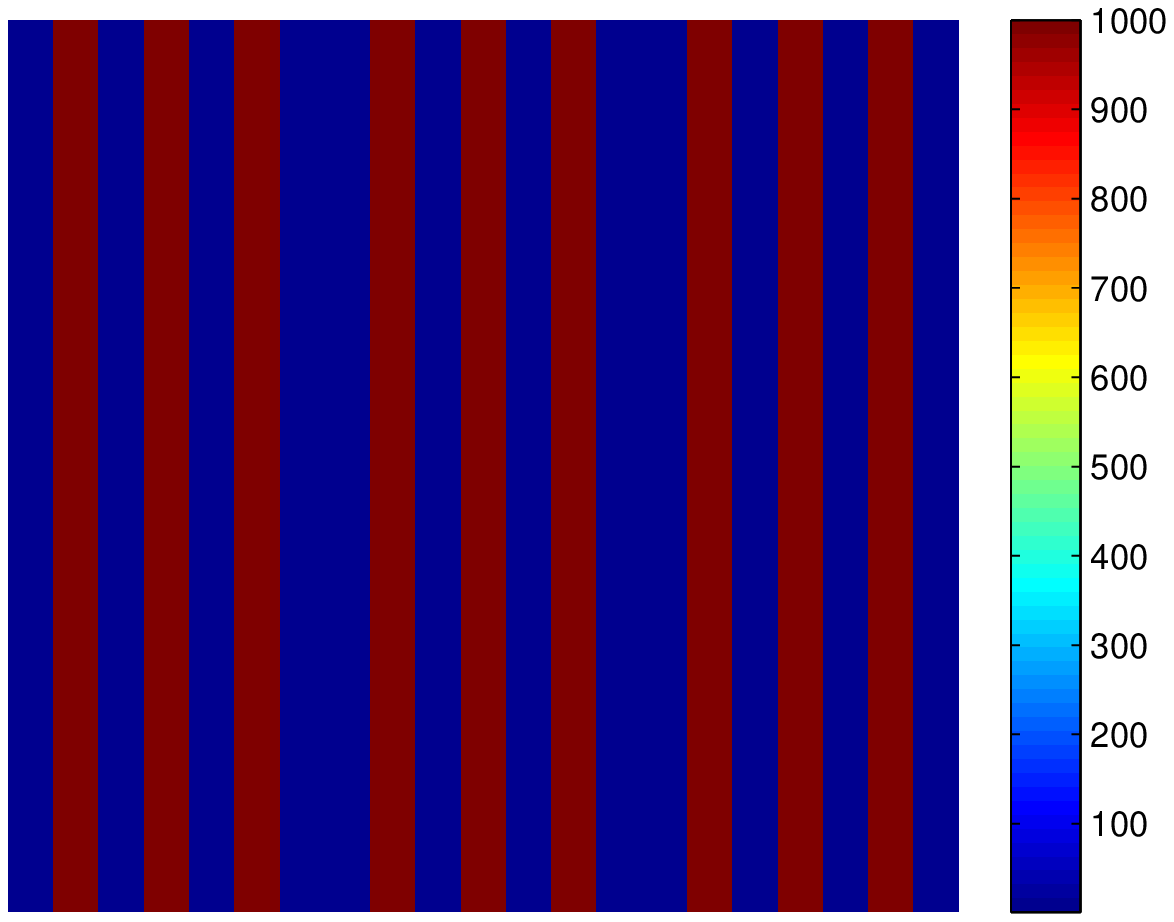}
\end{minipage}
    \caption{$N = 3^2$ with one channels (left) and three channels (right) in each subdomain: blue ($\rho(x) = 1$) and red ($\rho(x) = \eta$).}
  \label{fig-channels}
\end{figure}

\begin{example}\label{example-2}
Consider model problem \eqref{continuous-var} with $\rho(x)$, which has channel patterns as shown in Figure \ref{fig-channels}.
\end{example}

We fixed $N = 3^2$, $\beta=1/2$ and $\eta = 10^3$ , and the results for Example \ref{example-2} are presented in Table \ref{table-2-A}.
We note that the total number of primal unknowns in M1 increase as more channels are introduced,
 but it is nearly independent of the local problem size in both methods.
In particular, M2 chooses lesser primal unknowns than M1.

\begin{figure}[H]
   \centering
\begin{minipage}[t]{0.45\linewidth}
    \includegraphics[width=0.9\textwidth]{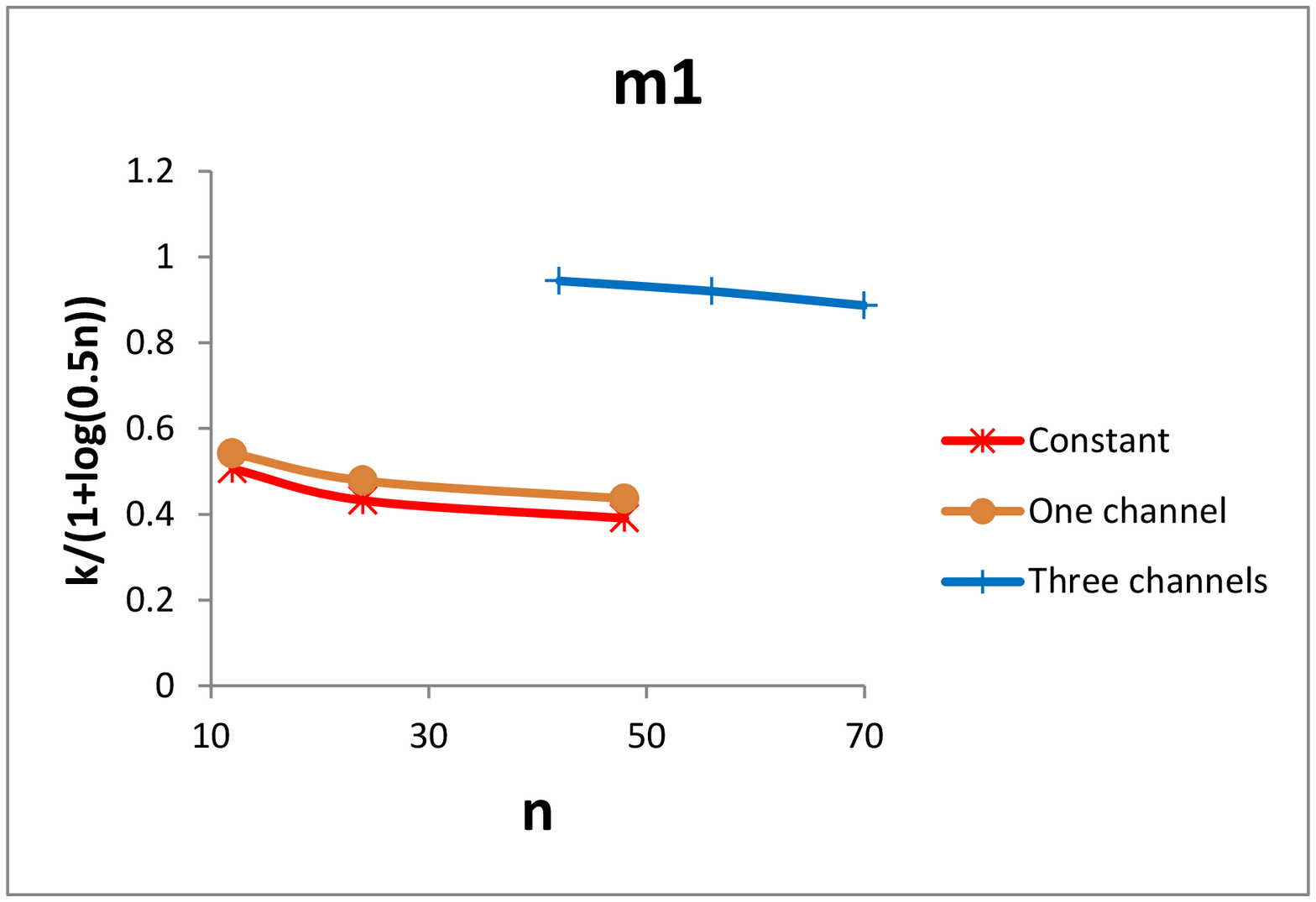}
\end{minipage}
\begin{minipage}[t]{0.45\linewidth}
    \includegraphics[width=0.9\textwidth]{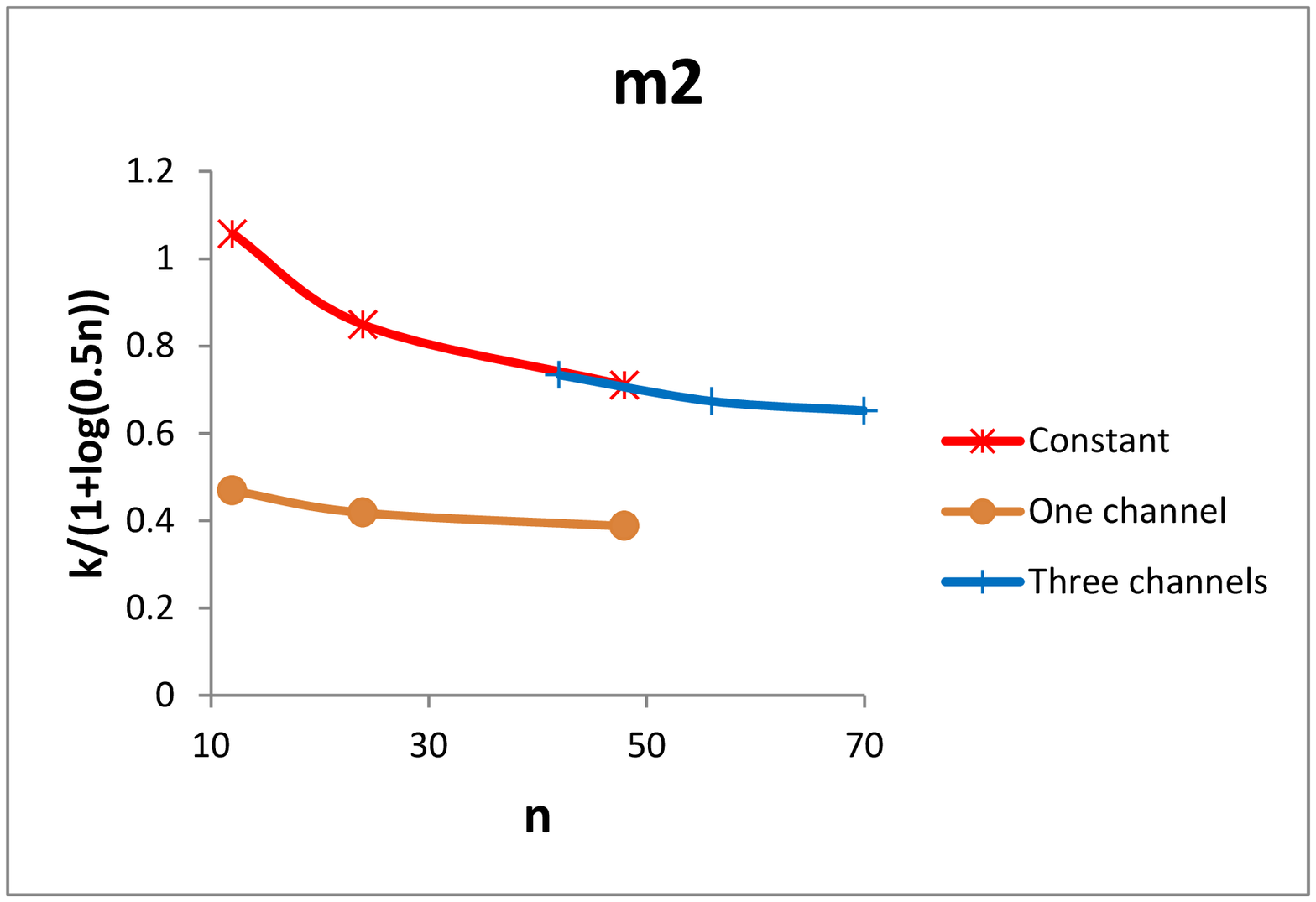}
\end{minipage}
    \caption{$\kappa/(1+log(0.5n))$ of M1 and M2 with varying $n$ for constant and channel $\rho(x)$}
  \label{fig-C}
\end{figure}

\begin{table}[H]
\centering\caption{Performance for Example \ref{example-2} with fixed $N = 3^2$, $\beta=1/2$ and $\eta = 10^3$}
\label{table-2-A}\vskip 0.1cm
\begin{tabular}{{|c|c|c|c|c|c|c|}}\hline
Channel &n                  &method   &Iter  & $\lambda_{\min}$   & $\lambda_{\max}$ &pnum  \\\hline
                      & \multirow{2}{*}{12}  &M1  &9   &1.0000   &1.4122   & 34\\
                      &                      &M2  &9   &1.0001   &2.9506   & 16\\\cline{2-7}
\multirow{3}{*}{one} &\multirow{2}{*}{24}    &M1  &8   &1.0001   &1.5060   & 34\\
                      &                      &M2  &9   &1.0000   &2.9579   & 16\\\cline{2-7}
                      & \multirow{2}{*}{48}  &M1  &9   &1.0001   &1.6317   & 34\\
                      &                      &M2  &9   &1.0001   &2.9668   & 16\\\hline\hline
                      & \multirow{2}{*}{42}  &M1  &10  &1.0001   &3.8197   & 66\\
                      &                      &M2  &11  &1.0000   &2.9666   & 16\\\cline{2-7}
\multirow{3}{*}{three} & \multirow{2}{*}{56} &M1  &11  &1.0001   &3.9869   & 64\\
                      &                      &M2  &11  &1.0000   &2.9183   & 16\\\cline{2-7}
                      & \multirow{2}{*}{70}  &M1  &11  &1.0001   &4.0407   & 64\\
                      &                      &M2  &11  &1.0000   &2.9701   & 16\\\hline
\end{tabular}
\end{table}
In Figure \ref{fig-C}, we plot $C = \kappa/\Theta$ of M1 and M2 with varying $n$ for the constant and channel $\rho(x)$, where $\Theta = 1+ log(0.5n)$. It is easy to see that the constant $C$ in Theorem \ref{condition-number} is independent of $n$.

For the case with three channels and $n=42$, we present the numerical results for varying $\eta$ in Table \ref{table-3-A}.
We can see that the two methods are both robust to $\eta$. Especially, as $\eta$ increases, the number of primal unknowns of M2 stays the same.

\begin{table}[H]
\centering\caption{Performance for Example \ref{example-2} with varying $\eta$ in three channels, and fixed $n=42$, $N=3$, $\beta = \frac{1}{2}$}
\label{table-3-A}\vskip 0.1cm
\begin{tabular}{{|c|c|c|c|c|c|c|c|}}\hline
$\eta$               &method   &Iter  & $\lambda_{\min}$   & $\lambda_{\max}$ &pnum \\\hline
\multirow{2}{*}{$10$}  &M1  &11   &1.0008   &1.9610 &16 \\
                       &M2  &9    &1.0001   &1.9325 &16 \\\hline
\multirow{2}{*}{$10^2$}&M1  &15   &1.0001   &3.9311 &34 \\
                       &M2  &10   &1.0000   &2.7299 &16 \\\hline
\multirow{2}{*}{$10^3$}&M1  &10   &1.0001   &3.8197 &66 \\
                       &M2  &11   &1.0000   &2.9666 &16 \\\hline
\multirow{2}{*}{$10^4$}&M1  &9    &1.0001   &1.5811 &70 \\
                       &M2  &11   &1.0000   &2.9953 &16 \\\hline
\multirow{2}{*}{$10^5$}&M1  &9    &1.0001   &1.6025 &70 \\
                       &M2  &12   &1.0000   &3.0008 &16 \\\hline
\end{tabular}
\end{table}

\begin{figure}[H]
  \centering
  \includegraphics[width=0.4\textwidth]{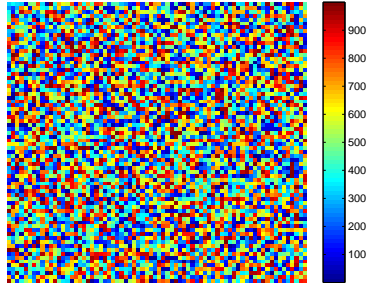}\\
  \caption{The coefficient for random $\rho(x)$ from $10^{-3}$ to $10^3$ with $N = 4^2$, $n = 18$ and $\beta = 1$.}\label{fig_3D_Element}
\end{figure}

\begin{example}\label{example-3}
Consider model problem \eqref{continuous-var} with $\rho(x) = 10^r$, where $r$ is chosen randomly from $(-3,3)$ for each grid element, as shown in Figure \ref{fig_3D_Element}.
\end{example}

For a given $\beta = \frac{3}{2}$, we present the numerical results of both methods for increasing $n$ with a fixed $N=3^2$ in Table \ref{table-4-A}, where the average number of primal unknowns per interface is given in the parentheses.
For M1, the number of adaptive primal unknowns is more than $50\%$ of the total interface unknowns.
But it's worth pointing out that M2 has a significant advantage in iteration number and gives about 2 primal unknowns per interface as $n$ increases,
which shows that M2 is more robust and efficient for highly random coefficients than M1.
In Table \ref{table-5}, the two methods are tested for highly varying and random $\rho(x)$ by increasing $N$ with a fixed $n = 24$. We observe a similar performance to the previous case.

\begin{table}[H]
\centering\caption{Performance for Example \ref{example-3} by increasing $n$ with a fixed $N = 3$ and $\beta = \frac{3}{2}$}
\label{table-4-A}\vskip 0.1cm
\begin{tabular}{{|c|c|c|c|c|c|c|}}\hline
n                  &method   &Iter  & $\lambda_{\min}$   & $\lambda_{\max}$  &pnum   & ppnum \\\hline
\multirow{2}{*}{12}  &M1  &19   &1.0000   &3.3817  & 183(15.25) &66.30\%\\
                     &M2  &12   &1.0003   &2.0596  &  18(1.50)  &6.52\%\\\hline
\multirow{2}{*}{24}  &M1  &22   &1.0001   &4.1523  & 371(30.92) &65.78\%\\
                     &M2  &14   &1.0008   &3.0392  &  21(1.75)  &3.72\%\\\hline
\multirow{2}{*}{48}  &M1  &24   &1.0003   &4.8344  & 650(54.17) &57.02\%\\
                     &M2  &15   &1.0009   &3.2978  &  19(1.58)  &1.67\%\\\hline
\end{tabular}
\end{table}

\begin{table}[H]
\centering\caption{Performance for Example \ref{example-3} by increasing $N$ with a fixed $n = 24$ and $\beta = \frac{3}{2}$}
\label{table-5}\vskip 0.1cm
\begin{tabular}{{|c|c|c|c|c|c|c|}}\hline
N                  &method   &Iter  & $\lambda_{\min}$   & $\lambda_{\max}$ &pnum & ppnum \\\hline
\multirow{2}{*}{$4^2$} &M1  &23   &1.0001   &4.1516    & 703(29.29)  &62.32\%\\
                       &M2  &16   &1.0008   &3.1044    &  48(2.00)   &4.26\%\\\hline
\multirow{2}{*}{$5^2$} &M1  &22   &1.0001   &4.1453    &1190(29.75)  &63.30\%\\
                       &M2  &17   &1.0004   &3.1094    &  86(2.15)   &4.57\%\\\hline
\multirow{2}{*}{$6^2$} &M1  &22   &1.0001   &4.1702    &1829(30.48)  &64.86\%\\
                       &M2  &19   &1.0005   &3.9451    & 136(2.27)   &4.82\%\\\hline
\end{tabular}
\end{table}

Then, the similar results of Example \ref{example-3} are also presented for the geometrically unconforming partitions.

\begin{figure}[H]
 \centering
\begin{minipage}[t]{0.49\linewidth}
    \includegraphics[width=0.85\textwidth]{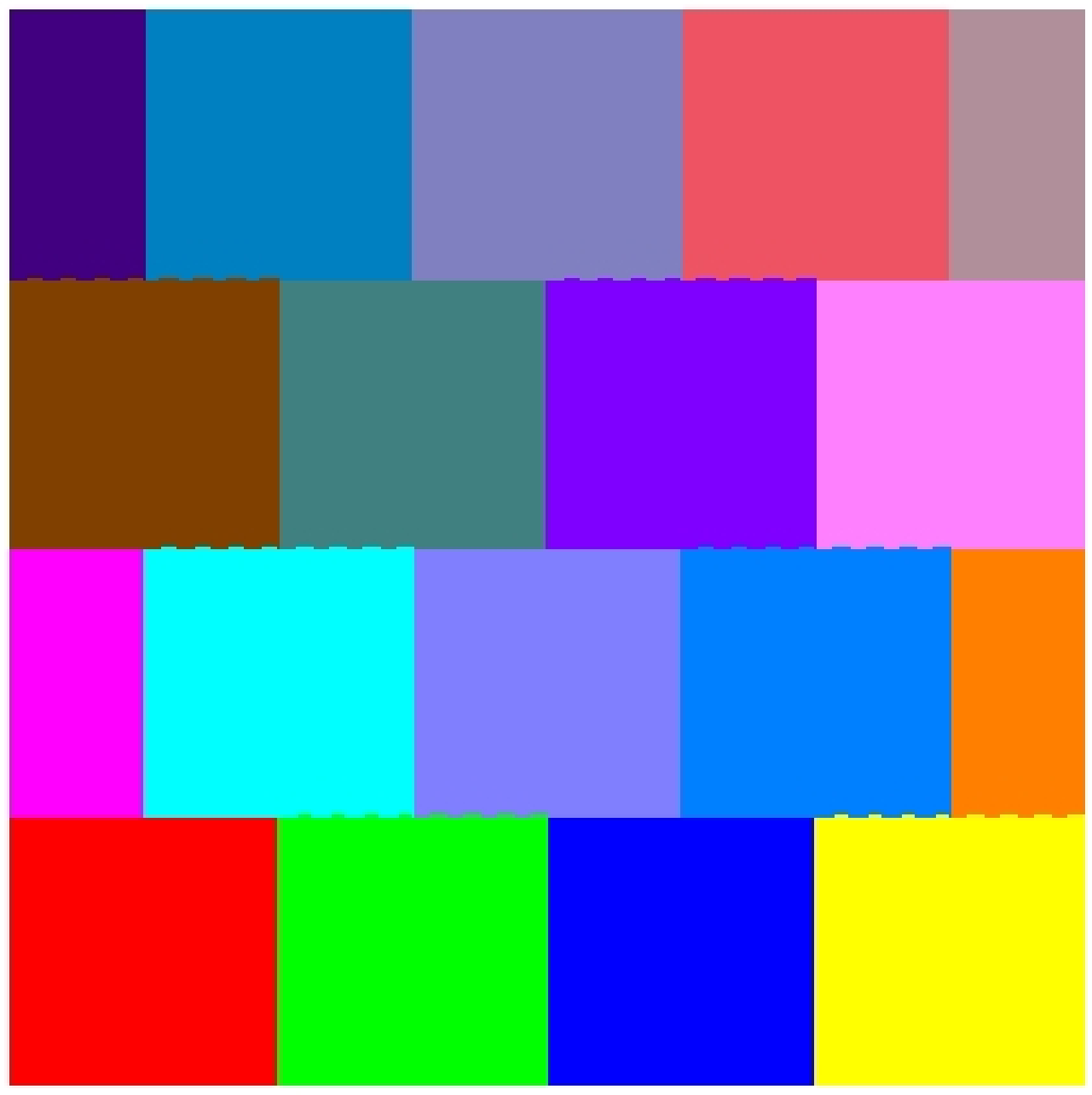}
\end{minipage}
\begin{minipage}[t]{0.49\linewidth}
    \includegraphics[width=0.75\textwidth]{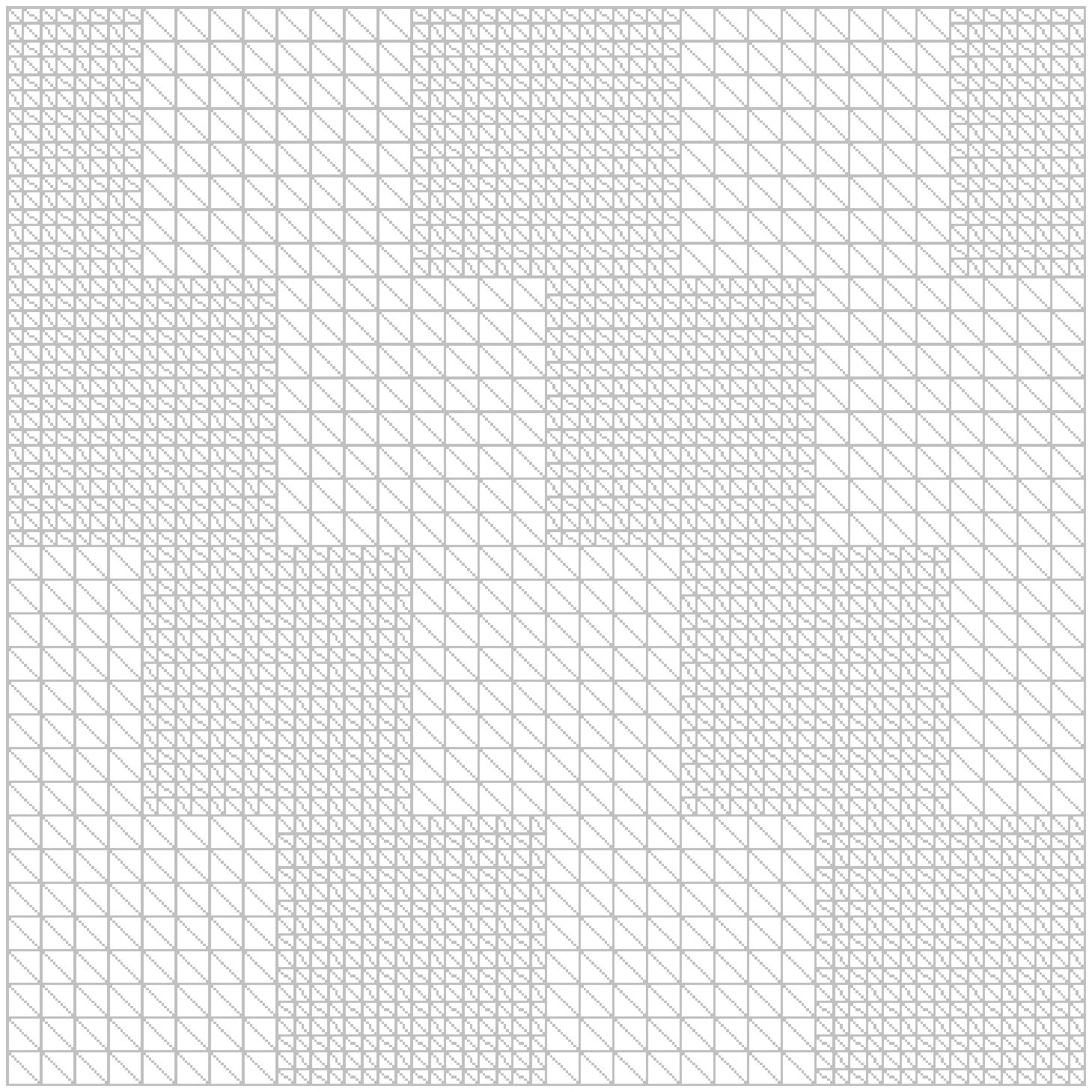}
\end{minipage}
    \caption{A geometrically unconforming partitions, $N = 18$, $n=8$ and $\beta = 2$.}
  \label{fig-unconform-111}
\end{figure}

For a given geometrically unconforming square partitions with $N=18$, see Figure \ref{fig-unconform-111},
without loss of generality, we assume that the space $X_h$ is associated with the $\mathcal{P}_1$ Lagrange finite element space,
the results for Example \ref{example-3} by increasing $n$ with $\beta = 2$ are shown in Table \ref{table-6}.

\begin{table}[H]
\centering\caption{Performance of the methods with a geometrically unconforming subdomain partition}
\label{table-6}\vskip 0.1cm
\begin{tabular}{{|c|c|c|c|c|c|c|}}\hline
n                  &method   &Iter  & $\lambda_{\min}$   & $\lambda_{\max}$ &pnum & ppnum \\\hline
\multirow{2}{*}{8} &M1  &16   &1.0000   &2.6274   & 131(3.45)   &67.53\%\\
                   &M2  &15   &1.0004   &2.2196   &  49(1.29)   &25.26\%\\\hline
\multirow{2}{*}{16}&M1  &20   &1.0000   &3.4079   & 233(6.13)   &54.69\%\\
                   &M2  &17   &1.0004   &3.1529   &  59(1.55)   &13.85\%\\\hline
\multirow{2}{*}{32}&M1  &22   &1.0001   &4.1848   & 390(10.26)  &43.82\%\\
                   &M2  &19   &1.0003   &3.5826   &  60(1.58)   &6.74\%\\\hline
\multirow{2}{*}{64}&M1  &25   &1.0002   &4.9182   & 698(18.37)  &38.39\%\\
                   &M2  &18   &1.0006   &3.3565   &  57(1.50)   &3.14\%\\\hline
\end{tabular}
\end{table}

\begin{remark}
For $\varepsilon = 0$, we can get the similar results by using the regularized techniques in \cite{HU2007}.
\end{remark}

From all the experiment results above, we find that the condition numbers confirm our theoretical estimate.
The two methods are all robust for the constant and channel $\rho(x)$, but for highly varying and random coefficients,
M2 shows better performance than M1.

\section{Conclusions}

In this paper, we develop an adaptive BDDC algorithm in variational form for
high-order mortar discretizations by introducing some vector-valued auxiliary spaces and operators with essential properties.
Since there is not any continuity constraints at subdomain vertices in the mortar method involved in this paper,
it simplifies the construction of the primal unknowns.
We show that the condition number of the preconditioned system is bounded by
a given tolerance, which is used to construct the transformation operators for selecting coarse basis functions.
Numerical results are presented to verify the robustness and efficiency of the proposed approaches.

\section*{Acknowledgements}
This work is supported by the National Natural Science Foundation of China
(Grant Nos. 11571293, 11201398, 11301448, 11601462),
Hunan Provincial Natural Science Foundation of China (Grant No. 2016JJ2129).

%

\bibliographystyle{elsarticle-num}

\FloatBarrier

\end{document}